\documentclass[regno, reqno,10pt]{amsart}
\usepackage{amsmath,amsthm,amssymb,amsfonts,epsfig,color,graphicx,enumerate,accents,subfigure,comment, esint}
\usepackage[a4paper,margin=2.59truecm]{geometry}
\usepackage{graphicx}
\usepackage{enumitem}

\usepackage[linkcolor=blue,colorlinks=true, citecolor = blue]{hyperref}

\DeclareMathOperator{\dive}{div}

\DeclareMathOperator{\diam}{diam}

\DeclareMathOperator{\dist}{dist}
\DeclareMathOperator{\supp}{supp}
\DeclareMathOperator{\Tr}{Tr}

\def\eps{{\varepsilon}}
\def\wto{\rightharpoonup}

\def\N{\mathbb{N}}

\def\R{\mathbb{R}}
\def\O{\Omega}
\def\o{\omega}

\def\vf{\varphi}

\def\HH{\mathcal{H}}

\def\FF{\mathcal{F}}
\def\GG{\mathcal G}

\def \< {\langle}
\def \>{\rangle}
\def\HH{\mathcal{H}}
\def\II{\mathcal{I}}
\def\JJ{\mathcal{J}}

\def\MM{\mathcal{M}}

\def\RR{\mathcal{R}}
\def\TT{\mathcal{T}}
\def\JJ{\mathcal{J}}
\def\Leb{\text{Leb}}
\def\f{\phi}

\newcommand{\be}{\begin{equation}}
\newcommand{\ee}{\end{equation}}

\newcommand{\bs}{\begin{split}}
\newcommand{\es}{\end{split}}

\newcommand{\red}{\textcolor{red}}

\newcommand{\LN}[2]{\| #1 \|_{L^N(#2)} }
\newcommand{\Linfty}[2]{\| #1 \|_{L^{\infty}(#2)} }
\renewcommand{\L}[3]{\| #2 \|_{L^{#1}(#3)}}
\renewcommand{\b}[1]{\bar{#1}}

\numberwithin{equation}{section}
\theoremstyle{plain}
\newtheorem{theorem}{Theorem}[section]
\newtheorem{lemma}[theorem]{Lemma}
\newtheorem{corollary}[theorem]{Corollary}
\newtheorem{proposition}[theorem]{Proposition}
\newtheorem{definition}[theorem]{Definition}
\theoremstyle{remark}
\newtheorem{remark}[theorem]{Remark}

\newenvironment{ack}{{\bf Acknowledgements. }}

\title[Optimal Regularity in Transmission Problems]{Optimal regularity for variational solutions of free transmission problems.}

\author{Diego Moreira$^{*}$}
\address{$^*$Departamento de Matemática, Universidade Federal da Ceara(Fortaleza, Brazil)}
\email{$^*$dmoreira@mat.ufc.br}
\author{Harish Shrivastava$^\dagger$}
\address{$^\dagger$Tata Institute of Fundamental Researcher-Centre of of Applicable Mathematics}
\email{$^\dagger$harish21@tifrbng.res.in, }

\begin{document}
\maketitle
\begin{abstract}
In this article we study functionals of the type considered in \cite{HS21}, i.e.
$$
J(v):=\int_{B_1} A(x,u)|\nabla u|^2 +f(x,u)u+ Q(x)\lambda (u)\,dx
$$
here $A(x,u)= A_+(x)\chi_{\{u>0\}}+A_-(x) \chi _{\{u<0\}}$, $f(x,u)= f_+(x)\chi_{\{u>0\}}+f_-(x) \chi _{\{u<0\}}$ and $\lambda(x,u) = \lambda_+(x) \chi_{\{u>0\}} + \lambda_-(x) \chi_{\{u\le 0\}}$.  We prove the optimal  $C^{0,1^-}$ regularity of minimizers of the functional indicated above (with precise H\"older estimates) when the coefficients $A_{\pm}$ are continuous functions and $\mu \le A_{\pm}\le \frac{1}{\mu}$ for some $0<\mu<1$, with $f \in L^N(B_1)$ and $Q$  bounded. We do this by presenting a new compactness argument and approximation theory similar to the one  developed by L. Caffarelli in \cite{Ca89} to treat the regularity theory for solutions to fully nonlinear PDEs. Moreover, we introduce the $\TT_{a,b}$ operator that allows one to transfer minimizers from the transmission problems to the Alt-Caffarelli-Friedman type functionals, {in small scales,} allowing this way the study of the regularity theory of minimizers of Bernoulli type free transmission problems. 

\end{abstract}

\medskip

\indent \textbf{2010 Mathematics Subjects Classification:} 35B65, 35Q35, 35R35, 35J20.

\tableofcontents
\section{Introduction}
As nicely pointed out in \cite{CSS21},  M. Picone introduced transmission problems in theory of elasticity in 1954 (\cite{picone}) and the theory was further developed by Lions \cite{L56}, Stampaccia \cite{S56} and Campanato \cite{C57}. The non-divergence case was treated by Schechter in 1960 (\cite{S60}). For further details, we refer to \cite{CSS21}, from where we learn about the history of those developments.   
\begin{comment}
{ \cite{LU68}, \cite{CF12, LN03, LV00, MOV09}.
\end{comment}

Mathematical models involving composite materials have gained significant attention in recent years as we can see, for instance, from the works of \cite{CSS21}, \cite{D21}, \cite{FKS20}, \cite{CKS21}, \cite{PS20}, \cite{Z21} (see also a slightly older and interesting article \cite{CF12}). It is a challenging task to model the electromagnetic or heat conduction properties of certain materials, particularly when they are close to threshold points which lead to abrupt changes in the nature of the material. The same issue arises when one tries to model properties of composite mixtures. Applications and results that deal with such class of models can be found in \cite{M10}, \cite{AI96} and others. Interesting discussions related to those models in the applied sciences are also present in \cite{BDHS81}, \cite{WLSXMZ14}, \cite{BC84}. 

In a very rough way, the models discussed above can be represented by the PDE
\be\label{PDE formulation}
\dive A(x,u,\nabla u) = f(x,u)
\ee
where $A(x,u,\nabla u)$ and $f(x,u)$ have jump discontinuities with respect to the variable $u$. These jump discontinuities model the change in the diffusion coefficients as soon as $u$ crosses a threshold value or else, the phase of the material changes. Few examples of phenomenas that can be modelled through transmission problems are mixture of different chemicals, conductivity (electric and thermal) of composite materials or a material operating close to thresholds like triple point, melting point or breakdown potential. 

Transmission problems can also be posed in variational formulations where configurations of least energy are of interest.  This means that one is led to study minimizers of functionals of the following type
\be\label{variational formulation}
\int_{\O}H(x,u,\nabla u)\,dx
\ee
where $H(x,u,\nabla u)$ has a jump discontinuity with respect to its variables. Some results related to the above mentioned variational formulation \eqref{variational formulation} can be  found in \cite{HS21}, \cite{TA15}. Moreover, the form \eqref{PDE formulation} is also studied in \cite{KLS17}, \cite{KLS19}.

In this work, {we study the optimal regularity result  for variational solutions to transmission problems of Bernoulli type}. In fact, in our context here, we have a free transmission problem in the sense that there is no a priori knowledge on position, geometry or regularity of the separating interface. More precisely, we focus on (local) minimizers of functionals of the form 
\be\label{P}
J_{A,f,Q} (v;B_1) =\int_{B_1} A(x,v)|\nabla v|^2+f(x,v)v+ Q(x)\lambda (v)\,dx
\ee
where for $v\in H^1(B_1)$ we set
\[
\begin{split}
A(x,v)&:= A_+(x)\chi_{\{v>0\}}+A_-(x) \chi _{\{v<0\}},\\
f(x,v)&:= f_+(x)\chi_{\{v>0\}}+f_-(x) \chi _{\{v<0\}},\\
\lambda(v) &:= \lambda_+ \chi_{\{v>0\}} + \lambda_- \chi_{\{v\le 0\}}.
\end{split}
\]
In a way, our transmission problem is a inhomogeneous version of the classical two phase free boundary problem studied by H. Alt, L. Caffarelli and A. Friedman in 1984 (c.f. \cite{ACF84}) where now the leading gradient term on the functional has jump discontinuities across the free boundary and also a inhomogeneity term (which is not present in \cite{ACF84}). In our context here, the following assumptions are enforced throughout the paper
\begin{enumerate}[label=\textbf{(G\arabic*)}]
\item\label{G1} (\textbf{Continuity}) $A_{\pm}$ are  continuous in $B_1$, the open unit ball in $\R^N$. 
\item\label{G2} $f_{\pm}\in L^N(B_1)$.
\item\label{G3} (\textbf{Ellipticity}) There exists $\mu\in (0,1)$ such that 
$$
\mbox{for every $x\in B_1$ we have $\mu  \le A_{\pm}(x) \le  \frac{1}{\mu}$.}
$$
\item \label{G4}  $\Lambda:= \lambda_+ - \lambda_->0$.
\item \label{G5} $0\le Q(x) \le q_2<\infty$ for every $x\in B_1$. Without loss of generality, we assume $q_2\ge 1$.
\end{enumerate}
\vspace{0.5cm}

Here, unlike \cite{ACF84} (also \cite{KLS17} for the non-variational setting), local minimizers are not locally Lipschitz continuous. In fact, we show here that they are $C^{0,1^-}_{loc}$ (with precise estimates) and this is the optimal regularity. As a matter of fact, this optimal regularity is already manifested at the PDE level as it is proven in \cite{JMS09} (see Remark \ref{low-reg} below).

\vspace{0.1cm}

Our strategy in this paper can be divided in two macro parts. In the first one, we obtain the optimal growth rate for local minimizers in balls centered on the zero level set. In the second part, we prove the optimal regularity of minimizers by using the Euler-Lagrange equations that are satisfied in both phases. We then match those regularities by using Harnack type arguments to obtain the final optimal regularity (see for instance, Proposition \ref{step 3.5}).
\vspace{0.1cm}

In order to implement the second part pointed out above, we also need analogous estimates for weak solutions of PDEs of the form \eqref{PDE} $($c.f. Theorem \ref{PDE regularity}$)$. We were unable to find such results in the literature showing how precisely the constants depend on the data of the problem as well as results that hold in any dimension $N\geq 2$. This is done in Section \ref{PDE section} of the paper. 

\vspace{0.1cm}
For the optimal growth rate of minimizers along the zero level set, we implement an approximation theory similar to the one developed by L. Caffarelli in the seminal paper \cite{Ca89} to treat regularity theory for solutions to fully non-linear PDEs. There are essentially two main steps. The first one is to prove the proximity of local minimizers to regular profiles by compactness argument (c.f Theorem \ref{step 1} and Proposition \ref{step 1.5}). In our case, the regular profiles inherit regularity from (local) minimizers of functionals of Alt-Caffarelli-Friedman type (as in \cite{ACF84}) via $\TT_{a,b}$ operator discussed in the sequel. The second one is to reduce the problem to a so called ``small regime configuration" (c.f. Proposition \ref{step 2}) via the scaling invariance of the problem under the appropriate regularity assumptions of the data.

We observe that compactness here appears to be subtler than the analogous situation for viscosity solutions in the theory PDEs, as stability properties of minimizers (by which we mean minimality property being kept under limits) are in general more delicate to assure.  Moreover, in the context of nonlinear variational Bernoulli free boundary type problems, local Lipschitz regularity (unavailable in our case here) appears to be an effective auxiliary instrument to prove stability (at least in the level of blow-up of minimizers) as one can see in the works of D. Danielli and A. Petrosyan (Lemma 5.3 in \cite{DP06}), S. Martinez and N. Wolanski (Lemma 7.2 in \cite{MW08}) and H.  Alt, L. Caffarelli, A. Friedman themselves (Lemma 3.3 in \cite {ACF84}).

\begin{comment}
Furthermore, it seems to us that in the context of nonlinear variational problems involving Bernoulli type free boundary problems, local Lipschitz regularity (unavailable in our case here) appears to be a recurrent idea to prove stability for blow-ups, especially in the non-linear case. This is the case for instance in the works of D. Danielli and A. Petrosyan (Lemma 5.3 in \cite{DP06}), S. Martinez and N. Wolanski (Lemma 7.2 in \cite{MW08}) and H.  Alt, L. Caffarelli, A. Friedman themselves (Lemma 3.3 in \cite {ACF84}). 
\end{comment}

Here, we follow a different route inspired on the ideas due to H. Alt and L. Caffarelli present in \cite[Lemma 5.4]{AC81} (for blow-ups of minimizers) that account for careful perturbations of minimizers to create suitable admissible competitors for energy comparison. Furthermore, these ideas (in \cite{AC81}) also account for smart cancelations driven by the weak convergence, showing somehow a dependence on the quadratic structure of the functional.

In our case, although we do not deal with blow-ups (of minimizers), we still need a careful construction in the spirit of the proof of \cite[Lemma 5.4]{AC81} in order to produce suitable competitors for minimizers. The construction is more involved since our problem is a nonlinear one (it is a transmission problem instead of a pure free boundary one of Bernoulli type). Moreover, the competitors need to be constructed in such a way that the $L^1_{loc}$ convergence of positivity sets holds, recovering this way some kind of  upper semi-continuity of the phase transition part of the functional, namely,  $\int Q(x)\lambda(v)\,dx$ \footnote{This is reflected, for instance, in equations \eqref{charv p=2} and \eqref{charv2} in the proof of Theorem \ref{step 1}.}. We close the compactness argument using the Widman's hole filling technique (Proposition \ref{step 1.5}).

Finally, we highlight that to handle the transmission problem \textit{per se}, we introduce a new tool, namely, the $\TT_{a,b}$ operator. The idea is that this operator ``bridges" minimizers of transmission functionals of the form \eqref{J0} with minimizers of functionals of Alt-Caffarelli-Friedman type in \cite{ACF84} (c.f Lemma \ref{corollary to lemma 1.3}). This permits a perturbation theory to be carried out (like in \cite{Ca89}) by passing from one functional to the other in small scales, recovering this way the regularity theory for the original minimizers.  We choose to present the results pertaining to the $\TT_{a,b}$ operator in Section \ref{Tab operator} for general value of $p\ge 1$, i.e, transmission functionals involving $p-$energy, as it can be of independent interest and of use in other transmission problems. At this point, we do not know how to implement the compactness argument for the case where $p\in [1,\infty)$ and $p\neq 2$ since the cancellation effect yielded by the weak convergence (see for instance, the proof of \ref{E4} in Theorem \ref{step 1}) seems to be delicate to reproduce in the case $p\neq 2$.

Our results extend the ones in \cite{TA15} and \cite{HS21} in the particular case where matrix coefficients $A_{\pm}$ are variable multiples of the identity matrix. Indeed, in those papers the authors prove that minimizers become more and more regular as the ``jump" distance between the phase coefficients $A_+$ and $A_-$ tends to zero. Here, as pointed out before, we show that any minimizer of \eqref{P} is locally $C^{0,1^-}$ regular with precise corresponding H\"older estimates, independent of the ``jump'' distance between $A_+$ and $A_-$ across the phases. This paper is essentially self contained and it has detailed proofs.

\vspace{0.1cm}
\textbf{Outline of the paper: }
In  Section \ref{setting the problem}, we list the main definitions and results. In Section \ref{Tab operator} we introduce the $\TT_{a,b}$ operator which converts a functional of the form \eqref{J0} into a functional of Alt-Caffrelli-Friedman type. We show that the $\TT_{a,b}$ operator preserves the regularity of minimizers of the original functional with constant coefficients $($c.f. Proposition \ref{lemma 1.3}$)$. Then in Section \ref{PDE section}, we prove the interior $C^{0,1^-}$ regularity estimates for weak solutions to the PDE of the form \eqref{PDE} $($c.f. Theorem \ref{PDE regularity}$)$. In Section \ref{compactness}, we make use of ideas from \cite{AC81} and show that under the small regime, minimizers of \eqref{P} are close to minimizers of functional of the form \eqref{J0}, whose regularity estimates are well understood thanks to Section \ref{Tab operator}. Finally, in Section \ref{main result} and \ref{main theorem} we prove the main result via rescaling arguments, namely Theorem \ref{C01-}.

\section{Main definitions and results}\label{setting the problem}

In this paper, we study the regularity of minimizers of the energy functional \eqref{P}, under the assumptions \ref{G1}-\ref{G5}. For the coming sections, we simplify the notation as needed by dropping some subscripts whenever no ambiguity is present. For instance, $J_{A,f,Q}(\cdot, B_1)$ can be denoted only by $J$ when definitions of $A,f,Q$ and $B_1$ are clear from the context. We will note $F$ as 
\be\label{def of F}
F:= |f_+| +|f_-|.
\ee
Let $D\subset \R^N$ be a bounded, Lipschitz and open set. Then for all $\f \in W^{1,p}(D)$ 
\be\label{bdry data space}
\begin{split}
W_{\phi}^{1,p}(D)\,&:=\, \left \{ v\in W^{1,p}(D)\,:\,\; v-\phi \in W_0^{1,p}(D) \right \}\\
&=\Big \{ v\in W^{1,p}(D)\,:\, \Tr_{\partial D}(v)=\Tr_{\partial D}(\f) \Big \}.
\end{split}
\ee
When $p=2$, $W^{1,2}(D)$ is denoted by $H^1(D)$ and analogously we define $H^1_{\f}(D)$. We say that a function $v\in W^{1,p}(D)$ takes the boundary value $\f$ on $\partial D$ whenever $v\in W_{\phi}^{1,p}(D)$. In order to study the traces, we use the notation $L^{p}(\partial D)$ for $L^{p}(\partial D,d\HH^{N-1})$.

Existence of a minimizer for $J_{A,f,Q}(\cdot, B_1)$ in \eqref{P} follows from \cite[Theorem 2.1]{HS21}. For the analysis that we carry out in later sections, we define the following notation (\ref{J0}) for $x_0\in D\subset B_1$
\be\label{J0}
J_{x_0}(v;D) :=\int_{D}A_+(x_0)|\nabla v^+|^2+A_-(x_0)|\nabla v^-|^2+Q(x_0)\lambda(v)\,dx.
\ee

\begin{definition}
Given a functional $J(\cdot;\O):W^{1,p}(\O)\to \R$, a function $u\in W^{1,p}(\O)$ is an absolute minimizer of $J (\cdot;\O)$ if for every $v\in W^{1,p}(\O)$ such that $u-v\in W_0^{1,p}(\O)$ we have 
$$
J (u;\O)\le J (v;\O).
$$
A function $u\in W^{1,p}(\O)$ is a local minimizer of $J (\cdot;\O)$ if  for every open Lipschitz subset $D\subset \subset \O$ and for every $v\in W^{1,p}(D)$ such that $\Tr_{\partial D}(v)=u$ we have 
$$
J (u;D)\le J (v;D).
$$
\end{definition}
In fact, absolute minimizers of functionals of the form \eqref{P} and \eqref{J0} are also local minimizers (c.f. Lemma \ref{local minimizer}). In order to better state our results, we define the modulus of continuity of continuous function $A:\O\to \R$ as follows

\begin{definition}\label{mod of cont}
A modulus of continuity is a function $\o:[\,0,\infty)\to \R$ such that $\o$ is non-negative, non-decreasing, and $\lim_{t\to 0} \o(t)=0$. In particular if $A: \O \to \R$ is continuous and $D\subset \subset \O$. We define $\o_{A,D}$ as 
$$
\o_{A,D}(t) :=\begin{cases} \sup \Big ( |A(x)-A(y)| \,:\,\; x,y\in D \;\;\mbox{such that $|x-y|\le t$} \Big )\;\;\;\;\;\;\;t\le  \diam(D)\\
\o_{A,D}(\diam(D))\qquad\qquad \qquad\qquad\qquad\qquad \qquad\qquad \;\;\;\;\; \;\;\;t> \diam (D).
\end{cases}
$$
Then $\o_{A,D}$ is a modulus of continuity $($we call $\o_{A,D}$ as the modulus of continuity of $A$ in $D$$)$.
\end{definition}
{We denote $C(\O)^{N\times N}$ as the set of $N\times N$ matrices with entries as continuous functions on $\O$.} 
\vspace{0.5cm}

\noindent \textbf{Main results:} We now state the main results in this paper.

\begin{theorem}[Optimal regularity for absolute minimizers] \label{C01-}
Let $u \in H^1(B_1)$ be a bounded absolute minimizer of $J_{A,f,Q}(\cdot, B_1)$ given in \eqref{P} and satisfying \ref{G1}-\ref{G5}. Then $u\in C_{loc}^{0,1^-}(B_1)$ with estimates. More precisely, for any $0<\alpha<1$ and $B_r \subset \subset B_1,$ we have that $u\in C^{0,\alpha}(B_r)$ with the following estimate

\be \label{main 1}
 \|u \|_{C^{\alpha}(B_r)} \le \frac{C(N, \alpha, \mu, q_2,\lambda_+, \o_{A_{\pm},B_{r^*}})}{(1-r)^\alpha}\left ( 1+ \Linfty{u}{B_1}+\LN{f}{B_1} \right ). 
 \ee
Here $r^*= \frac{1+r}{2}$. The function $\o_{A_{\pm}, B_{r^*}}$ is the maximum of modulus of continuity of $A_{+}$ and $A_-$ in $B_{r^*}$. 
\end{theorem}

In the process of proving the main result, we obtain optimal PDE regularity result which seems to us to be unavailable in the literature in that form (where there is a precise constant dependence and the result holding in  every dimension $N\ge2$)
\begin{theorem}[Optimal regularity of bounded weak solutions]\label{PDE regularity}
Let $u\in H^{1}(B_1)$ be a bounded weak solution to
$$
\dive\left (A(x)\nabla u\right )=f\;\;\mbox{in $B_1$}
$$
where $A\in C(B_1)^{N\times N}$ and there exists $0<\mu<1$ such that $\mu |\xi|^2 \le \< A(x)\xi ,\xi \>\le \frac{1}{\mu}|\xi|^2$ a.e. in $B_1$ for all $\xi \in \R^N$, $f\in L^N(B_1)$. Then $u\in L^{\infty}_{loc}(B_1)$ and for any $0<\alpha<1$ and $B_r \subset \subset B_1,$ we have that $u\in C^{0,\alpha}(B_r)$ with the following estimate
\be \label{main 0.1}
 \|u \|_{C^{\alpha}(B_r)} \le \frac{C(N, \alpha, \mu, \o_{A,B_{r^*}})}{(1-r)^\alpha}\left ( \Linfty{u}{B_1}+\LN{f}{B_1} \right ). 
 \ee
Here $r^*= \frac{1+r}{2}$. The function $\o_{A, B_{r^*}}$ is the modulus of continuity of $A$ in $B_{r^*}$. 

\end{theorem}

\begin{remark}\label{low-reg}
The regularity provided in Theorem \ref{PDE regularity} and Theorem \ref{C01-} are sharp even in the case where right hand side $f$ is equal to zero. For this, we refer to \cite{JMS09}. For recent developments in the subject, we refer to \cite{MM21}.
\end{remark}
\begin{comment}
Next, we present a scaled version of our result. For that, we define the following defintion 
\begin{definition}
Let $v\in C^{0,\alpha}(B_r)$, then we define $\|\cdot \|^{*}_{C^{0,\alpha}(B_r)}$ as
$$
\|v\|^*_{C^{0,\alpha}(B_r)} = \Linfty{v}{B_r}+ r^{\alpha}[v]_{C^{0,\alpha}(B_r)}.
$$ 
\end{definition}
\begin{corollary}
Let $u\in H^{1}(\O)$ be a local minimizer of $J_{A,f,Q}(\cdot, \O)$. Then $u$ is locally bounded in $\O$ and for every $x_0\in \O$ we have 
$$
\|u \|^*_{C^{\alpha}(B_{r}(x_0))} \le {C(N,p, \alpha, \mu, q_2\lambda_+,\o_{A_{\pm}, B_{2r}(x_0)})}\left ( r+\Linfty{u}{B_{2r}(x_0)}+r \LN{f}{B_{2r}(x_0)} \right ).
$$
Here $r<\frac{d}{4}$ $(d:= \dist(x_0,\partial \O))$, $\o_{A_{\pm}, B_{2r}(x_0)}$ is the modulus of continuity of $A_{\pm}$ in the ball $B_{2r}(x_0)$. 
\end{corollary}
\end{comment}
\section{Preliminary tool: The $\TT_{a,b}$ operator}\label{Tab operator}
In this section $D$ is a bounded, open and Lipschitz subset on $\R^N$. We start by considering the following functionals
$$
\JJ(v;D):= \int_{D} \Big ( a^p|\nabla v^+|^p+b^p|\nabla v^-|^p+\lambda(v) \Big ) \,dx \mbox{ and }\FF(v;D):=\int_{D} \Big ( |\nabla v|^p+\lambda(v) \Big )\,dx,
$$
where $p\in (1,\infty)$ and $a,b>0$. We make the observation that when $u$ is a minimizer  of the energy functional $\JJ(\cdot; ;D)$, then, this very function carry the regularity properties available for a minimizer of the classical functional $\FF(\cdot ;D)$. The purpose of this section is to develop the ingredients that prove the above remark. For that, we introduce the $\TT_{a,b}$-operator.
\begin{definition}\label{Tab}
Let $a,b>0$, $p\in [1,\infty)$ and $D$ be an open Lipschitz set, we define $\TT_{a,b}:W^{1,p}(D)\to W^{1,p}(D)$ as follows
$$
\TT_{a,b}(v) := a v^+\,-\,bv^-.
$$
\end{definition}
\begin{remark}[Definition of $\TT_{a,b}$-operator on $\partial D$]\label{Tab trace}
In the coming sections, we deal with the action of the $\TT_{a,b}$-operator on functions in $W^{1,p}(D)$ and also keep track of their traces. For this reason, in order to maintain compatibility of notations, we define the $\TT_{a,b}$-operator also on the boundary level acting on $L^p(\partial D)$ to be  $\TT^{\partial D}_{a,b}:L^p(\partial D) \to L^p(\partial D)$ given by 
$$
\TT^{\partial D}_{a,b}(\varphi) := a \varphi^+\,-\,b\varphi^-.
$$
\end{remark}
%In order to simplify the notation, we use simply denote $\TT^{\partial D}_{a,b}$ by $\TT_{a,b}$, whenever convinient.  
\begin{remark}\label{subsequence limits}
In the sequel, we use the following simple fact from the basic theory of topological spaces. 
Consider $(X,\tau)$ to be a Hausdorff topological space and let $v_k$ be a sequence in $X$ and $v\in X$. Suppose that for every subsequence $v_{k_j}$, we can extract yet another subsequence $v_{k_{j_l}}$ so that $v_{k_{j_l}} \to v$ in $(X,\tau)$ as $l\to \infty$. Then $v_k\to v$ in $(X,\tau)$ as $k\to \infty$.
\end{remark}

\begin{remark}\label{uniqueness of decomp}
We recall the following decomposition for any $a\in \R$. There exists unique numbers $p_a$ and $n_a$ such that $a=p_a-n_a$ and $p_a,n_a\ge 0$ and $p_a\cdot n_a =0$. Clearly, $p_a= \max \left ( a,0 \right )=:a^+$ and $n_a = -\min \left ( a,0 \right ) =:a^-$. 

This decomposition also applies to functions. For $w\in L^p(D)$, $\TT_{a,b}(v) = [\TT_{a,b}(v)]^+ - [\TT_{a,b}(v)]^-$ a.e. in $D$ and from the above mentioned uniqueness of decomosition
\be\label{unique pm decomp}
\begin{split}
[\TT_{a,b}(v)]^+& =  av^+=(av)^+\mbox{ a.e. in $D$}\\
[\TT_{a,b}(v)]^- &= b v^-= (bv)^-\mbox{ a.e. in $D$}.
\end{split}
\ee
Similarly, \eqref{unique pm decomp} hold for every $w\in L^p(\partial D)$ $\HH^{N-1}$ a.e. on $\partial D$.  Moreover,
$$
\chi_{\{v>0\}} = \chi_{\{\TT_{a,b}(v)>0\}}\,\mbox{ and }\chi_{\{v\le0\}} = \chi_{\{\TT_{a,b}(v)\le0\}} \mbox{ a.e. in $D$}.
$$
In particular, we have 
\be\label{lambda preserved}
\lambda(v) = \lambda(\TT_{a,b}(v)) \;\mbox{a.e. in $D$}.
\ee
\end{remark}
We also note the following observation
\begin{remark}\label{plus minus}
For every $w\in W^{1,p}(D)$, we observe that 
\be\label{abs val}
\begin{split}
|w^+\pm w^-|^p &= |w^+\pm w^-|^p\chi_{\{w>0\}} + |w^+\pm w^-|^p\chi_{\{w\le 0\}}\\
&= \big |(w^+\pm w^-)\chi_{\{w>0\}}\big |^p + \big |(w^+\pm w^-)\chi_{\{w\le 0\}}\big |^p\\
&= |w^+|^p + |w^-|^p.
\end{split}
\ee
Also,
\be\label{abs val nabla}
\begin{split}
|\nabla w^+\pm \nabla w^-|^p &= |\nabla w^+\pm \nabla w^-|^p\chi_{\{w>0\}} + |\nabla w^+\pm \nabla w^-|^p\chi_{\{w\le 0\}}\\
&= \big |(\nabla w^+\pm \nabla w^-)\chi_{\{w>0\}}\big |^p + \big |(\nabla w^+\pm \nabla w^-)\chi_{\{w\le 0\}}\big |^p\\
&= |\nabla w^+|^p + |\nabla w^-|^p.
\end{split}
\ee
The above two computations imply that 
\be\label{Lp plus minus}
\begin{split}
\L{p}{w}{D}^p=\L{p}{ w^+\pm  w^-}{D}^p = \L{p}{w^+}{D}+ \L{p}{w^-}{D}^p\\
\end{split}
\ee
\be\label{Lp nabla plus minus}
\begin{split}
\L{p}{\nabla w}{D}^p=\L{p}{ \nabla w^+\pm \nabla  w^-}{D}^p = \L{p}{\nabla w^+}{D}+ \L{p}{\nabla w^-}{D}^p
\end{split}
\ee
\end{remark}
The above remarks lead to the following important corollary
\begin{corollary}\label{Lp Tab}
Let $p\in (1,\infty)$ then,  the $\TT_{a,b}$ operator splits as follows
\begin{enumerate}[label = $($\alph*$)$]
\item \label{a} For every $w\in L^{p}(D)$ 
\be\label{Tab Lp}
\begin{split}
|\TT_{a,b}w|^p &= a^p|w^+|^p+b^p|w^-|^p\mbox{ a.e. in $D$}\\
\end{split}
\ee
in particular
$$
\L{p}{\TT_{a,b}(w)}{D}^p = a^p \L{p}{w^+}{D}^p + b^p \L{p}{w^-}{D}^p.
$$
\item\label{b} \eqref{Tab Lp} also holds for $\vf\in L^p(\partial D)$ replacing the Lebesgue measure by $H^{N-1}_{\lfloor_ {\partial D}}$.
\vspace{0.2cm}
\item \label{c} For every $w\in W^{1,p}(D)$
\be\label{Tab Lp nabla}
\begin{split}
|\nabla \TT_{a,b}(w)|^p &= a^p|\nabla w^+|^p + b^p|\nabla w^-|^p\mbox{ a.e. in $D$}\\
\end{split}
\ee
in particular
$$
\L{p}{\nabla \TT_{a,b}(w)}{D}= a^p \L{p}{\nabla w^+}{D}^p + b^p \L{p}{\nabla w^-}{D}^p.
$$
\item\label{d} $w\in W^{1,p}(\O)$ if and only if $\TT_{a,b}(w)\in W^{1,p}(\O)$.
\end{enumerate}
\end{corollary}
\begin{proof}
The claims \ref{a}, \ref{b} and \ref{c} follow readily by applying \eqref{abs val}, \eqref{abs val nabla}, \eqref{Lp plus minus} and \eqref{Lp nabla plus minus} for $\TT_{a,b}(w)$ and then using \eqref{unique pm decomp}. \ref{d} readily follows from \ref{a} and \ref{c}.
\end{proof}
We show that the $\TT_{a,b}$ operator is sequentially continuous in weak and strong $W^{1,p}(D)$ topologies and commutes with the trace operator. Moreover, it preserves the pointwise convergence of sequence of functions. The proofs of these facts are inspired by the {ideas from \cite{MM78} and \cite{JK}.}
\begin{proposition}\label{prop of Tab}The following properties hold 

\begin{enumerate}[label = \textbf{(\arabic*)}]

\vspace{0.2cm}

\item \label{P2} Assume $v_k \wto v$ $($weakly$)$ in $W^{1,p}(D)$. Then

\vspace{0.2cm}

\begin{enumerate}[label=\textbf{(1\alph*)}]
\item\label{p2}  $v_k^{\pm} \wto v^{\pm} \,\mbox{$($weakly$)$ in $W^{1,p}(D)$}. $
\vspace{0.3cm}

\item \label{p2.1} $\TT_{a,b} (v_k) \wto \TT_{a,b} (v)$ $($weakly$)$ in $W^{1,p}(D)$. 
\vspace{0.3cm }
\item \label{p2.2} The maps $v\mapsto |v|$ and $\TT_{a,b}$ operator are sequentially weakly continuous in $W^{1,p}(D)$.
\vspace{0.1cm}
\end{enumerate}

\vspace{0.4cm}

\item \label{P1} Assume $v_k\to v$ $($strongly$)$ in $W^{1,p}(D)$. Then 
\vspace{0.2cm}
\begin{enumerate}[label=\textbf{(2\alph*)}]
\item \label{p1}$v_k^{\pm} \to v^{\pm} \,\mbox{$($strongly$)$ in $W^{1,p}(D)$}.$
\vspace{0.3cm}
\item \label{p1.1} $\TT_{a,b}(v_k) \to \TT_{a,b}(v) $ $($strongly$)$ in $W^{1,p}(D)$. 
\vspace{0.3cm}
\item \label{p1.2} The maps $v\mapsto |v|$ and $\TT_{a,b}$ operator are continuous in $W^{1,p}(D)$.

\end{enumerate}

\vspace{0.2cm}

\item \label{p0} Assume $v_k\to v$ in $L^p(D)$. Then 
\vspace{0.2cm}
\begin{enumerate}[label=\textbf{(3\alph*)}]
\item\label{3a} $v_k^{\pm}\to v^{\pm} \;\mbox{in $L^p(D)$}.$
\vspace{0.3cm}
\item\label{3b} $\TT_{a,b} (v_k) \to \TT_{a,b}(v) \;\mbox{in $L^p(D)$}.$
\end{enumerate}
\vspace{.2cm}
\item \label{P3} Assume $v\in W^{1,p}(D)$. Then 
\vspace{0.2cm}
\begin{enumerate}[label=\textbf{(4\alph*)}]
\item \label{p4} $\Tr(v^{\pm}) = \left ( \Tr(v) \right )^{\pm}\,\mbox{ $\HH^{N-1}$ a.e. on $\partial D$}.$
\vspace{0.3cm}
\item \label{p5}$\Tr(\TT_{a,b}(v)) = \TT^{\partial D}_{a,b} (\Tr(v)) \,\mbox{ $\HH^{N-1}$ a.e. on $\partial D$}.$

\end{enumerate}

\vspace{0.4cm}
\item \label{p6} \ref{p2}, \ref{p2.1}, \ref{p1} and \ref{p1.1} also hold when $W^{1,p}(D)$ is replaced by $W^{1,p}_0(D)$.

\vspace{0.4cm}

\item \label{p7} Assume $v_k \to v$ pointwise almost everywhere in $D$. Then $\TT_{a,b}(v_k) \to \TT_{a,b}(v)$ pointwise almost everywhere in $D$. Moreover, if $v\in C(D)$, then $|v|, \TT_{a,b}(v) \in C(D)$.
\end{enumerate}
\end{proposition}

\begin{remark}
Although in \ref{p2.2} we showed the sequential weakly continuity of $\TT_{a,b}$-operator in $W^{1,p}(D)$, we in fact, can prove that $\TT_{a,b}$-operator is weakly continuous. For that we need to repeat the same arguments replacing the sequences by nets.
\end{remark}

\begin{proof}

In order to prove \ref{P2} we observe that from the assumptions that $v_k\in W^{1,p}(D)$ is a bounded sequence in $W^{1,p}(D)$. Moreover,   
\[
\begin{split}
\int_{D}|{v_k^{+}}|^p\,dx\,&\le \, \int_{D}|v_k|^p\,dx\\
\int_{D}|\nabla {v_k^{+}}|^p\,dx &=\int_{D}|\nabla v_k|^p\chi_{\{v_k>0\}}\,dx \,\le \, \int_{D}|\nabla v_k|^p\,dx.
\end{split}
\]
Hence, $v_k^{+}$ is also a bounded sequence in $W^{1,p}(D)$. Therefore there exists $w\in W^{1,p}(D)$ such that up to a subsequence which we still denote as $v_k^+$ 
\[
\begin{split}
v_{k}^{+} \wto w \,&\mbox{ weakly in $W^{1,p}(D)$ and }\\
v_k^+ \to w\,&\mbox{ in $L^p(D)$}.
\end{split}
\]
From compact embedding of $W^{1,p}(D)$ in $L^p(D)$, we know that $v_k\to v$ in $L^p(D)$. Furthermore, we observe that 
\be\label{claim 0.0}
\int_{D}|v_k^{\pm}-v^{\pm}|^p\,dx \le  \int_{D}|v_k-v|^p\,dx\to 0\;\mbox{as $k\to \infty$}.
\ee
From \eqref{claim 0.0}, we conclude $v_k^+\to v^+$ in $L^p(D)$. Therefore, from Remark \ref{subsequence limits} we have $v_k^+\to v^+=w$ in $L^p(D)$. So, 
$$
v_k^+\wto  v^+\,\mbox{ weakly in $W^{1,p}(D)$.}
$$
Similarly, we can prove that $v_k^-\wto  v^-\,\mbox{ weakly in $W^{1,p}(D)$.}$ This proves \ref{p2}, \ref{p2.1} and \ref{p2.2} follow readily.
In order to prove \ref{P1}, from \eqref{Lp nabla plus minus} and \cite[Lemma 7.6, Lemma 7.7]{GT98} we observe that  for every $w\in W^{1,p}(D)$,
\be\label{a1}
\begin{split}
\int_{D}\Big |\nabla \left (w^\pm \mp\frac{w}{2} \right )\Big |^p\,dx &= \int_{D} \Big | \nabla \left ( \frac{w^+}{2} + \frac{w^-}{2}\right ) \Big |^p\,dx\\
&= \left (\frac{1}{2}\right )^p \int_{D}|\nabla w|^p\,dx\\
\end{split}
\ee
Since $v_k\to v$ in $W^{1,p}(D)$, we have 
\be\label{a2}
\begin{split}
\lim_{k\to \infty} \int_{D} \Big | \nabla \left ( v_k^\pm \mp\frac{v_k}{2} \right ) \Big |^p\,dx &= \lim_{k\to \infty}\left ( \frac{1}{2} \right )^p \int_{D}|\nabla v_k|^p\,dx \\
&= \left ( \frac{1}{2} \right )^p \int_{D}|\nabla v|^p\,dx\\
&= \int_{D} \Big | \nabla \left ( v^\pm \mp\frac{v}{2} \right ) \Big |^p\,dx.
\end{split}
\ee
From \ref{p2} we know that 
$$
\left ( v_k^\pm \mp\frac{v_k}{2} \right ) \wto \left ( v^\pm \mp\frac{v}{2} \right )\;\;\mbox{weakly in $W^{1,p}(D)$}.
$$
From weak convergence and covergence of norms \eqref{a2} by  \cite[Chapter 6, Proposition 9.1, pp. 259]{EDB16} (Radon's theorem), we obtain
$$
\left ( v_k^\pm \mp\frac{v_k}{2} \right ) \to \left ( v^\pm \mp\frac{v}{2} \right )\mbox{ strongly in $W^{1,p}(D)$. }
$$
Therefore, $v_k^\pm \to v^\pm\;\;\mbox{in $W^{1,p}(D)$}.$
This proves \ref{p1}, \ref{p1.1} and \ref{p1.2} follow readily. Moreover, the argument in \eqref{claim 0.0} also proves \ref{3a} and also \ref{3b}.

Now, let us discuss the proof of \ref{p4} and \ref{p5}. Let $v\in W^{1,p}(D)$. From \cite[Theorem 4.3]{GE15} there  exists a sequence $v_k\in C(\overline D)\cap W^{1,p}(D)$ such that $v_k\to v$ in $W^{1,p}(D)$. Then from \cite[Theorem 4.6]{GE15}, we have 
\be\label{0}
\Tr(v_k^\pm) =  v_k^{\pm} \lvert _{\partial D} =  \Tr(v_k)^\pm \;\mbox{ on $\partial D$.}
\ee
From \ref{p1}, we know that $v_k^\pm \to v^\pm$ in $W^{1,p}(D)$ and then
\be\label{1}
\begin{split}
&\Tr(v_k^\pm) \to \Tr(v^{\pm})\mbox{  in $L^p(\partial D)$}.\\
\end{split}
\ee
\begin{comment}
\red{and therefore 
\be\label{1.1}
\Tr(v_k^\pm) \wto \Tr(v^{\pm})\mbox{ weakly  in $L^p(\partial D)$}.
\ee
}
\end{comment}
Also, by continuity of trace operator $\Tr(v_k)\to \Tr(v)$ in $L^p(\partial D)$. Now by \eqref{0} and from the argument  \eqref{claim 0.0}  with Lebesgue measure replaced by $d\HH^{N-1}$ restricted to $\partial D$ we obtain
\be\label{2}
\begin{split}
&\Tr(v_k)^\pm \to \Tr(v)^{\pm}\mbox{  in $L^p(\partial D)$}.\\
\end{split}
\ee
Combining the equations (\ref{0}), (\ref{1}) and (\ref{2}) we conclude 
$$
\Tr(v)^\pm  \stackrel{\text{\eqref{2}}}{=} \lim_{k\to \infty} \Tr(v_k)^{\pm} \stackrel{\text{\eqref{0}}}{=}  \lim_{k\to \infty} \Tr(v_k^{\pm}) \stackrel{\text{\eqref{1}}}{=} \Tr(v^{\pm})\mbox{  $\;\;\HH^{N-1}$ a.e. on $\partial D$}.
$$
Commutation of the Trace operator $\Tr$ and $\TT_{a,b}$-operator (i.e. \ref{p5})  follows directly from \ref{p4}.

In order to see \ref{p6} we note that since $v_k\in W_0^{1,p}(D)$, $\Tr(v_k)=0$ by \cite[Section 5.5, Theorem 2]{evans}. From \ref{p4} $\Tr(v_k^{\pm})=0$ and thus $v_k^{\pm}\in W_0^{1,p}(D)$ by \cite[Section 5.5, Theorem 2]{evans}. Now, \ref{p6} follows from the fact that $W_0^{1,p}(D)$ is closed and weakly closed in $W^{1,p}(D)$ \cite[Theorem 3.7]{HB10}.

Finally we have \ref{p7} by observing that 
$$
v_k\to v \mbox {  a.e. in $D$} \implies v_k^\pm \to v^\pm \mbox {  a.e. in $D$} \implies \TT_{a,b}(v_k) \to \TT_{a,b}(v) \mbox {  a.e. in $D$}.
$$

\end{proof}

\begin{proposition}\label{one one onto}
Regarding the operator $\TT_{a,b}:W^{1,p}(D)\to W^{1,p}(D)$. We have 
\vspace{0.3cm}
\begin{enumerate}[label = \textbf{(\arabic*)}]
\item \label{q1} $\TT_{a,b}: W^{1,p}(D) \to W^{1,p}(D)$ is a homeomorphism and 
$$
\Big [ \TT_{a,b} \Big ]^{-1} = \TT_{\frac{1}{a},\frac{1}{b}}.
$$
\item \label{q1.5} $\TT_{a,b}^{\partial D}: L^{p}(\partial D) \to L^{p}(\partial D)$ is a homeomorphism and 
$$
\Big [ \TT^{\partial D}_{a,b} \Big ]^{-1} = \TT^{\partial D}_{\frac{1}{a},\frac{1}{b}}.
$$
\item \label{q2} $\TT_{a,b}$ and $\TT_{a,b}^{-1}$ are bounded maps in $W^{1,p}(D)$ and $L^p(D)$ with the following estimates, 
\\[5pt]
\[
\begin{split}
&\min (a,b) \cdot \|w \|_{L^{p}(D)} \le \|\TT_{a,b}(w) \|_{L^{p}(D)}\le \max(a,b)\cdot \|w \|_{L^{p}(D)}\\[5pt]
& \min (a,b) \cdot \|\nabla w \|_{L^{p}(D)} \le \|\nabla \TT_{a,b}(w) \|_{L^{p}(D)}\le \max(a,b) \cdot\|\nabla w \|_{L^{p}(D)}\;\;\;\forall w\in W^{1,p}(D).
  \end{split}
\]

\hspace{0.8cm}

\noindent  Furthermore, $\TT_{a,b}^{\partial D}$ and $(\TT_{a,b}^{\partial D})^{-1}$  are bounded maps in $L^{p}(\partial D)$ with the following estimates, 
$$
  \min (a,b) \cdot \|w \|_{L^{p}(\partial D)} \le \|\TT^{\partial D}_{a,b}(w) \|_{L^{p}(\partial D)}\le \max(a,b) \cdot\|w \|_{L^{p}( \partial D)}\;\;\forall w\in L^p(\partial D).
$$
\item \label{q3} Let $\f\in W^{1,p}(D)$. Then, the image of $W^{1,p}_{\f}(D)$ under the operator $\TT_{a,b}$ is $W^{1,p}_{\TT_{a,b}(\f)}(D)$. Moreover, the map 
$$
\TT_{a,b}: W^{1,p}_{\f}(D)\to W^{1,p}_{\TT_{a,b}(\f)}(D)
$$
is a homeomorphism.
\end{enumerate}
\end{proposition}
\begin{proof}
In order to prove \ref{q1}, it is enough to show that 
\be\label{Tab invertible}
\TT_{\frac{1}{a},\frac{1}{b}} \left ( \TT_{a,b} \right ) =\TT_{a,b} \left (\TT_{\frac{1}{a},\frac{1}{b}}\right ) = Id \mbox{ in $W^{1,p}(D)$}.
\ee
Let $v\in W^{1,p}(D)$, then for any $\alpha ,\beta >0$
\be\label{Tab inj comp}
\begin{split}
\TT_{\frac{1}{\alpha},\frac{1}{\beta}} \Big ( \TT_{\alpha,\beta}(v) \Big )& = \frac{1}{\alpha} \Big [ \TT_{\alpha,\beta}(v) \Big ]^+ - \frac{1}{\beta} \Big [ \TT_{\alpha,\beta}(v) \Big ]^-\\
&=\frac{1}{\alpha} (\alpha v^+) -\frac{1}{\beta}(\beta v^-) =v.
\end{split}
\ee
The last identity follows from \eqref{unique pm decomp}. 
The claim \eqref{Tab invertible} readily follows by putting the pairs $(\alpha, \beta) = (a,b)$ and then $(\alpha, \beta)= (\frac{1}{a}, \frac{1}{b})$ in \eqref{Tab inj comp}.

Since the operators $\TT_{a,b}$ and $\Big [ \TT_{a,b} \Big ]^{-1} = \TT_{\frac{1}{a},\frac{1}{b}}$ are continuous in $W^{1,p}(D)$ (Proposition \ref{prop of Tab}-\ref{P1}) then $\TT_{a,b}$-operator is a homeomorphism. The proof of \ref{q1.5} follows similarly. The proof of \ref{q2} follows from Corollary \ref{Lp Tab}. 

In order to prove \ref{q3}, we take $\f\in W^{1,p}(D)$. We only need to show that 
\be\label{equal sets}
\TT_{a,b}\left (  W^{1,p}_{\f}(D)\right ) = W^{1,p}_{\TT_{a,b}(\f)}(D).
\ee
Let $\alpha, \beta>0$ and $w\in W^{1,p}_{\f}(D)$. Clearly, $\TT_{a,b}(w)\in W^{1,p}(D)$. Now, we compute the trace of $\TT_{a,b}(w)$. By Proposition \ref{prop of Tab}-\ref{p5}, we have
$$
\Tr(\TT_{\alpha,\beta}(w))= \TT^{\partial D}_{\alpha,\beta}(\Tr(w))= \TT_{\alpha,\beta}^{\partial D}(\Tr(\f)) =\Tr(\TT_{\alpha,\beta}(\f)).
$$ 
From this, we conclude that
\be\label{equal sets 1}
\TT_{\alpha,\beta}\left (  W^{1,p}_{\f}(D)\right ) \subset W^{1,p}_{\TT_{\alpha,\beta}(\f)}(D).
\ee
Therefore by taking the pair $(\alpha, \beta) = (a,b)$ we arrive to 
\be\label{blue star}
\TT_{a,b}\left (  W^{1,p}_{\f}(D)\right ) \subset W^{1,p}_{\TT_{a,b}(\f)}(D)
\ee
From Proposition \ref{one one onto}-\ref{q1} there exists a unique $\psi \in W^{1,p}(D)$ such that 
\be\label{phi inv}
\TT_{\frac{1}{a},\frac{1}{b}}(\psi)=\f.
\ee
Then plugging $(\alpha,\beta)= (\frac{1}{a},\frac{1}{b})$ and replacing $\f$ by $\psi$ in \eqref{equal sets 1}  
$$
\TT_{\frac{1}{a},\frac{1}{b}}\left (  W^{1,p}_{\psi}(D)\right ) \subset W^{1,p}_{\TT_{\frac{1}{a},\frac{1}{b}}(\psi)}(D).
$$
applying the $\TT_{a,b}$-operator on both sides of set inclusion above, using \eqref{phi inv} we arrive at
\be\label{equal sets 2}
\begin{split}
W^{1,p}_{\psi}(D)=\TT_{a,b} \left ( \TT_{\frac{1}{a},\frac{1}{b}}\left (  W^{1,p}_{\psi}(D)\right )  \right ) \subset \TT_{a,b} \left ( W^{1,p}_{\TT_{\frac{1}{a},\frac{1}{b}}(\psi)}(D) \right ) = \TT_{a,b} \left (  W^{1,p} _{\f}(D)\right ).
\end{split}
\ee
Recalling definition of $\psi$ in \eqref{phi inv} and from \eqref{equal sets 2} we obtain 
\be\label{equal sets 2.1}
\begin{split}
W^{1,p}_{\TT_{a,b}(\f)}(D) \subset \TT_{a,b} \left (  W^{1,p} _{\f}(D)\right ).
\end{split}
\ee
Finally, combining \eqref{blue star} and \eqref{equal sets 2.1} we obtain \eqref{equal sets}.
\end{proof}
%The following two propositions are the main result of this section
\begin{lemma}\label{corollary to lemma 1.3}
Consider the two functionals defined in $W^{1,p}(D)$, $p\in (\,1,\infty\,)$ 
$$
\FF(v;D):= \int_{D}|\nabla v|^p+\lambda(v)\,dx,\;\;\JJ(v;D):= \int_{D}a^p |\nabla v^+|^p+ b^p|\nabla v^-|^p +\lambda(v)\,dx.
$$
Then $u_0\in W^{1,p}(D)$ is an absolute minimizer of $\JJ(\cdot;D)$ if and only if $\TT_{a,b}(u_0)$ is an absolute minimizer of $\FF(\cdot;D)$.
\end{lemma}
\begin{proof}
Minimality of $\TT_{a,b}(u_0)$ for $\FF(\cdot;D)$ can be written as 
$$
\int_{D}|\nabla \TT_{a,b}(u_0)|^p+\lambda(\TT_{a,b}(u_0))\,dx\le \int_{D}|\nabla v|^p+\lambda(v)\,dx
$$
for every $v\in W^{1,p}(D)$ such that $v-\TT_{a,b}(u_0)\in W_0^{1,p}(D)$ or $v\in W^{1,p}_{\TT_{a,b}(u_0)}(D)$. From Proposition \ref{one one onto}-\ref{q3} and Proposition \ref{prop of Tab}-\ref{P3}, if $w$ runs over $W_{u_0}^{1,p}(D)$ if and only if $\TT_{a,b}(w)$ runs over $W_{\TT_{a,b}(u_0)}^{1,p}(D)$. Replacing $v$ by $\TT_{a,b}(w)$ in the inequality above, we obtain  
$$
\int_{D}|\nabla \TT_{a,b}(u_0)|^p+\lambda(\TT_{a,b}(u_0))\,dx\le \int_{D}|\nabla \TT_{a,b}(w)|^p+\lambda(\TT_{a,b}(w))\,dx
$$
for every $w\in W_{u_0}^{1,p}(D)$. Now, using \eqref{Tab Lp nabla} and recalling the equation \eqref{lambda preserved} $\lambda(\TT_{a,b}(u_0)) = \lambda(u_0)$ and $\lambda(\TT_{a,b}(w)) = \lambda(w)$, the proof of Proposition \ref{corollary to lemma 1.3} is finished.
\end{proof}

We call the functionals of the type \eqref{P}, \eqref{J0} or $\FF$ and $\JJ$ defined in the beginning of this section as "integral functionals". More precisely, for any $j: \O \times \R \times \R^N \to \R$ the funtional $J$ for the following form 
$$
J(v;\O):= \int_{\O}j(x,v,\nabla v)\,dx
$$
is called an integral functional. In order to proceed with the proof of the next lemma, we require the integral functional $J(\cdot;D)$ and  $J(\cdot;\O\setminus D)$ to be finite for every compact set $D\subset \subset \O$ for every $v\in W^{1,p}(\O)$.

\begin{lemma}\label{local minimizer}
 Assume $\O$ be an open, bounded and Lipschitz subset in $\R^N$ and $u\in W^{1,p}(\O)$ is an absolute minimizer of an integral functional $J(\cdot;\O)$ such that $J(v;D)<\infty$ for every $v\in W^{1,p}(\O)$ and $D\subset\subset  \O$. Then, $u$ is a local minimizer of $J(\cdot;\O)$.
\end{lemma}
\begin{remark}\label{finite J}
Lemma \ref{local minimizer} applies to the functionals of the form \eqref{P}, \eqref{J0}.  Indeed, assume $\O$ be an open and bounded subset of $\R^N$. We can easily see that for any $D\subset \subset \O$, $x_0\in D$ and $v\in H^1(D)$ 
$$
\max \Big ( |J_{A,f,Q}(v; D)|, |J_{x_0}(v;D)|,  |J_{A,f,Q}(v; \O\setminus D)|, |J_{x_0}(v;\O\setminus D)| \Big )<\infty.
$$
For more precise details, we refer to the computations in \eqref{fpmv p=2} and \eqref{fpmv p=2 N=2}.
\end{remark}

\begin{proof}[Proof of Lemma \ref{local minimizer}]
In order to prove that $u$ is a local minimizer of $J(\cdot;\O)$, we take $D\subset \subset \O$ an open Lipschitz subset. Let $v\in W^{1,p}(D)$ such that $\Tr_{\partial D}(v)=u$. We now set $w$ as 
$$
w=\begin{cases}
v\;\;\mbox{in $D$}\\
u\;\;\mbox{in $\O\setminus D$}.
\end{cases}
$$
Since $D$ is a Lipschitz domain and $\Tr_{\partial D}(v)=u$, from \cite[Theorem 3.44]{DD12} $w\in W^{1,p}(\O)$.  $\big ($Although \cite[Theorem 3.44]{DD12} is stated for $C^1$ domains, the same proof applies to Lipschitz domains by using the trace theorem \cite[Theorem 4.6]{GE15}$\big )$. Moreover, by the definition of $w$ and local properties of traces \cite[Theorem 5.7]{GE15}, $\Tr_{\partial \O}(w)=u$. Thus from the assumption we have 
$$
J(u;\O)\le J(w;\O).
$$
Now,
$$
J(u;D)+J(u;\O \setminus D) \le J(w;D)+J(w;\O \setminus D).
$$
From the definition of $w$, $w=u$ in $\O\setminus D$ and from Remark \ref{finite J} we cancel the second terms on both sides.  We obtain
$$
J(u;D) \le J(v;D).
$$
This proves that $u$ is a local minimizer of $J(\cdot; \O)$.
\end{proof}

\begin{proposition}\label{lemma 1.3}
For $p\in (1,\infty)$, assume $u_0\in W^{1,p}(D)\cap L^{\infty}(D)$ is an absolute minimizer of $\JJ(\cdot;D)$ where 
$$
\JJ(v;D):= \int_{D}a^p |\nabla v^+|^p+ b^p|\nabla v^-|^p +\lambda(v)\,dx.
$$
Then, $u_0\in C_{loc}^{0,1^-}(D)$ with estimates. 
\\[2pt]
\indent More precisely, for every $0<\alpha<1$ and $D'\subset \subset D$ we have the following estimates
\be\label{estl}
\|u_0\|_{C^{0,\alpha}(D')}\le \frac{1}{\min (a,b)}[\TT_{a,b}(u_0)]_{C^{0,\alpha}(D)}\le \frac{1}{\min(a,b)}C(N,\alpha, \lambda_+,\dist (D',\partial D), \Linfty{u_0}{D}).
\ee
\end{proposition}

\begin{proof}
From Lemma \ref{corollary to lemma 1.3} we know that $\TT_{a,b}(u_0)$ is an absolute minimizer of the $\FF(\cdot;D)$ defined in Lemma \ref{corollary to lemma 1.3}. This means
$$
\int_{D} |\nabla (\TT_{a,b}(u_0)) |^2\,dx +\lambda(\TT_{a,b}(u_0))\,dx \le \int_{D} |\nabla v |^2\,dx +\lambda(v)\,dx \;\;\;\forall \,v\in W^{1,p}_{\TT_{a,b}(u_0)}(D).
$$
From Lemma \ref{local minimizer} $\TT_{a,b}(u_0)$ is a local minimizer of $\FF(\cdot;D)$ as in Lemma \ref{corollary to lemma 1.3}. Hence, from \cite[Theorem 5.1]{JD14} we know that $\TT_{a,b}(u_0)\in C_{loc}^{0,1^-}(D)$ with the following estimates. 
$$
[\TT_{a,b}(u_0)]_{C^{\alpha}(D')} \le C(N, \mu, \lambda_+, \alpha, \dist(D',\partial D), \Linfty{u_0}{D})
$$
$\big ($One phase version of \cite[Theorem 5.1]{JD14} can be found in \cite[Theorem 4.1]{MW08}$\big )$.
This implies that $u_0$ also belongs to $C_{loc}^{0,1^-}(D)$ with estimates
\be
[u_0]_{C^{\alpha}(D')} \le \frac{1}{\min (a,b)}C(N, \lambda_+, \alpha,\dist(D',\partial D),\Linfty{u_0}{D}).
\ee
Indeed, let $x,y\in D'\subset \subset D$. Then for any given $0<\alpha<1$, we have
\be\label{TabLip}
|\TT_{a,b}(u_0)(x)-\TT_{a,b}(u_0)(y)| \le [\TT_{a,b}(u_0)]_{C^{0,\alpha}(D)}|x-y|^{\alpha}.
\ee
From Proposition \ref{prop of Tab}-\ref{p7}, we know that $\TT_{a,b}$-operator preserves continuity. We rewrite the above inequality in different regions of the domain $D'$. It is easy to see that 
\begin{itemize}
\vspace{0.1cm}
\item $\bar x \in \{u_0\ge 0\},\, \bar y\in \{u_0\ge0\} \implies  |u_0(\bar x)-u_0(\bar y)|\le \frac{[\TT_{a,b}(u_0)]_{C^{0,\alpha}(D')}}{a} |\bar x-\bar y|^\alpha$.
\vspace{0.3cm}
\item $\bar x \in \{u_0\le 0\},\, \bar y\in \{u_0\le0\} \implies  |u_0(\bar x)-u_0(\bar y)|\le \frac{[\TT_{a,b}(u_0)]_{C^{0,\alpha}(D')}}{b} |\bar x-\bar y|^\alpha$.
\vspace{0.3cm}
\item $\bar x \in \{u_0\ge0\},\, \bar y\in \{u_0\le0\} \implies  |u_0(\bar x)-u_0(\bar y)|\le \frac{[\TT_{a,b}(u_0)]_{C^{0,\alpha}(D')}}{\min(a,b)} |\bar x-\bar y|^\alpha$.
\end{itemize}
\vspace{0.2cm}
To check the last inequality we observe that if $\bar x \in \{u_0>0\},\, \bar y\in \{u_0\le0\}$ then
\vspace{0.2cm}
\[
\begin{split}
\min (a,b) \big | u_0(\bar x) - u_0(\bar y)\big | & \le \big ( a u_0(\bar x) - b u_0(\bar y) \big )\,\mbox{ (since the term inside is non-negative and $u_0(\bar y)\le 0$)}\\[5pt]
&=|au_0^+(x)-b(-u_0^-(x))|\\[5pt]
&= \big | \TT_{a,b}u_0(\bar x) - \TT_{a,b}u(\bar y) \big | \le [\TT_{a,b}(u_0)]_{C^{0,\alpha}(D')}|\bar x- \bar y|^{\alpha}.
\end{split}
\]
Then we combine the above inequality with \eqref{TabLip}.
\end{proof}
We present a basic lemma useful in upcoming proofs. The weak convergence in Lebesgue spaces preserve pointwise almost everywhere inequality.

\begin{lemma}\label{weak convergence preserves inequality}
Let $U\subset\mathbb{R}^N$ be a Lebesgue measurable set of finite measure and $f_k$ and $g_k$ be two sequences of functions in $L^p(U)$ with $p\in [1,\infty)$ such that 
$$
f_k\wto f\mbox{ and } g_k\wto g\;\mbox{ weakly in $L^p(U)$} \ \textnormal { and } \ f_k \le g_k \mbox{ a.e. in $U$ for all $k\in \N$.} 
$$
Then 
$$
f\le g \mbox{ a.e. in $U$}.
$$
\end{lemma}
\begin{proof}
For any $E\subset U$ Lebesgue measurable, $\chi_{E} \in L^{p'}(U)$. Then by the weak convergence 
$$
\int_{D}f \chi_{E}\,dx=\lim_{k\to \infty} \int_{D} f_k \chi_{E}\,dx \le  \lim_{k\to \infty} \int_{D} g_k \chi_{E}\,dx = \int_{D} g \chi_{E}\,dx.
$$
Therefore, $\int_{E}f\,dx \le \int_{E}g\,dx$ for every Lebesgue measurable set $E\subset U$. Hence $f\le g$ a.e in $U$.
\end{proof}

\section{Regularity of solutions of PDEs with continuous coefficients}\label{PDE section}

In this section we set $A(x)\in C^{N\times N}(D)$ and $f\in L^N(B_1)$. We are concerned about the regularity of (weak) solutions to the following PDE in $B_1$
\be\label{PDE}
\dive(A(x)\nabla u)=f.
\ee
Our setup condition regarding the above equation in this paper are
\begin{itemize}
\item \textbf{(Ellipticity)} There exists $0<\mu<1$ such that for almost every $x\in B_1$ and for all $\xi \in \R^N$ we have 
$$
\frac{1}{\mu} |\xi|^2 \le \< A(x)\xi, \xi \> \le \mu |\xi|^2.
$$
\item \textbf{(Continuity)} $A_{i,j} \in C(B_1)$ for every $1\le i,j\le N$.
\item $f\in L^N(B_1)$.
\end{itemize}
We remind the reader on the definition of weak solution to \eqref{PDE}.
\begin{definition}
$u\in H^1(B_1)$ is a weak solution to \eqref{PDE} in $B_1$ if for all $\vf \in H_0^{1}(B_1)$ we have
\be\label{weak solution}
\int_{B_1} \Big ( A(x)\nabla u\cdot \nabla \vf \Big )\,dx = -\int_{B_1} f(x) \vf(x)\,dx
\ee
\end{definition}
\begin{remark}
Let $h\in H_{loc}^{1}(B_1)\cap L^{\infty}(B_{1/2})$ such that $h$ weak solution to the following PDE in $B_{1/2}$
\be\label{p-harmonic}
\dive(A(0)\nabla h) =0.
\ee
Then $h$ satisfies the following classical interior gradient estimate
\be\label{int est}
\Linfty{\nabla h}{B_{1/4}} \le C(N, \mu)\Linfty{h}{B_{1/2}}.
\ee
\end{remark}
\begin{remark}\label{rescaling 0.1}
Let the function $u\in H^1_{loc}(B_{\Theta}(x_0))$ be a weak solution to the PDE \eqref{PDE} in $B_{\Theta}(x_0)$, we set $w$ as follows
$$
w(y) := \Phi u(\Theta y +x_0) + \Psi,  \;\;y\in B_1.
$$
It is easy to check by scaling and change of variables, $w\in H^1_{loc}(B_1)$ and $w$ is a weak solution to the following PDE 
$$
\dive(\bar A(x)\nabla w) = \bar f \;\;\mbox{in $B_1$}
$$
where $\bar A$ and $\bar f$ are defined as follows for $x\in B_1$
\[
\begin{split}
\bar A(x) &:= A(x_0+\Theta x )\\
\bar f(x)&:= \Phi \Theta ^{2} f(x_0+\Theta x ).
\end{split}
\]
\end{remark}

\begin{proposition}\label{step 0.1}
Assume $u\in H^{1}(B_{1/2})\cap L^{\infty}(B_1)$ be a weak solution to \eqref{PDE} in $B_{1/2}$ with $\Linfty{u}{B_{1/2}}\le 1$, $\L{2}{\nabla u}{B_{1/2}}\le M$ and $f\in L^N{(B_{1/2})}$. Then for every $\eps>0$ there exists $\delta(\eps, M,N, \mu)>0$ such that if 
$$
\max\left ( \|A(x) - A(0)\|_{L^{\infty}(B_1)}, \LN{f}{B_{1/2}} \right )\le \delta
$$
then 
$$
\Linfty{u-h}{B_{1/4}}\le \eps,
$$
for some $h\in H^1(B_{1/2})$ such that $\Linfty{h}{B_{1/2}}\le 1$ and $h$ is a weak solution to $\dive\left (A(0)\nabla h \right ) =0 $ in $B_{1/2}$.
\end{proposition}

\begin{proof}
Let us suppose for the sake of contradiction that the statement of Proposition \ref{step 0.1} is not true. This implies that there exists $\eps_0>0$ and sequences $A_k, f_k$ such that $A_{k}\in C(B_{1/2})$ and $f_k\in L^{N}(B_{1/2})$ satisfying $\|A_{k}-A(0)\|_{L^{\infty}(B_{1/2})}\le\frac{1}{k}$, $\|f_{k}\|_{L^N(B_{1/2})}\le \frac{1}{k}$. So that the corresponding weak solutions $u_k$ of the following PDE 
\be\label{PDEk}
\dive(A_k(x)\nabla u_k) = f_k(x)\;\;\mbox{in $B_{1/2}$}
\ee
are such that $\L{2}{\nabla u_k}{B_{1/2}}\le M$ and $\| u_k \|_{L^{\infty}(B_{1/2})} \le 1$ and for any $h$ satisfying the following PDE
\be\label{limit PDE}
\dive\left (A(0)\nabla h \right ) =0
\ee
we have
\be\label{absurdum}
\|u_k-h\|_{L^{\infty}(B_{1/4})}>\eps_0.
\ee
{Now, we see that $u_k$ are uniformly equicontinuous since they satisfy a uniform interior H{\"older} estimates. Indeed, from \cite[Corollary 4.18]{HL11} we know that there exists $\beta_0:= \beta_0 ( \mu,N)\in (0,1)$ and $C_0:= C_0 (\mu, N)>0$ so that for every $k\in \N$, $u_k$ satisfy the following uniform estimate 
\be\label{uk holder est}
|u_k(x)-u_k(y)|\le C_0  |x-y|^{\beta_0}   \Big \{ \L{2}{u_k}{B_{1/2}} + \LN{f_k}{B_{1/2}} \Big \}, \quad \forall x,y\in B_{1/4}.
\ee
}
Since $\L{2}{u_k}{B_{1/2}}\le C(N)\cdot \Linfty{u_k}{B_{1/2}}\le C(N)$ and $\LN{f_k}{B_{1/2}}\le \frac{1}{k}\le 1$, we can write
$$
\| u_k\|_{C^{\beta_0}(B_{1/4})}  \le C_1(N,\mu).
$$
By Arzela Ascoli theorem, there exists $u_0\in C^{0,\beta_0}(B_{1/4})$ such that 
\be\label{uniform convergence p=2} 
u_k\to u_0 \mbox{ in $L^{\infty}(B_{1/4})$ up to a subsequence. }
\ee
Since $\L{2}{\nabla u_k}{B_{1/2}}\le M$ and $\L{2}{u_k}{B_{1/2}} \le C(N)\Linfty{u_k}{B_{1/2}}\le C(N)$ for every $k$, we conclude that $u_k$ is a bounded sequence in $H^{1}(B_{1/2})$.  That is to say
$$
\| u_k \|_{H^{1}(B_{1/2})} \le C(M, N).
$$
Therefore, $u_k$ converges weakly to $u_0\in H^1(B_{1/2})$ (up to a subsequence)
$$
u_k \wto u_0 \;\;\mbox{weakly in $H^{1}(B_{1/2})$}.
$$
We claim that $u_0$ is a weak solution to $\dive(A(0)\nabla u_0)=0$, this gives us a contradiction and hence proves the Proposition \ref{step 0.1}.  Indeed, we know that $\L{2}{\nabla u_0}{B_{1/2}}\le M$ from lower semicontinuity of $H^1$ norm and  $\Linfty{u_0}{B_{1/2}}\le 1$ from Lemma \ref{weak convergence preserves inequality}. Once $u_k$ are weak solutions of \eqref{PDEk}, we have for every $\Phi \in H_0^{1}(B_{1/2})$
\be\label{e1}
\int_{B_{1/2}} \big ( A_k(x) \nabla u_k\cdot \nabla \Phi \big )\,dx = -\int_{B_{1/2}}f_k \Phi \,dx
\ee 
Therefore
\be\label{e2.}
\int_{B_{1/2}} \Big (  \big ( A_k(x)-A(0) \big )\nabla u_k\cdot \nabla \Phi \Big )\,dx + \int_{B_{1/2}} \big (A(0 ) \nabla u_k\cdot \nabla \Phi \big ) \,dx  = -\int_{B_{1/2}}f_k \Phi \,dx
\ee
Since $A_{k, \pm} \to A_{\pm}(0)$ uniformly in $B_{1/2}$ we obtain 
\[
\begin{split}
\int_{B_{1/2}} \Big (  \big ( A_k(x)-A(0) \big ) \nabla u_k\cdot \nabla \Phi \Big )\,dx & \le \Linfty{A_k-A(0)}{B_{1/2}} \int_{B_{1/2}}\big | \nabla u_k\cdot \nabla \Phi \big |\,dx\\
&\le \frac{1}{k} \,\,\L{2}{\nabla u_k}{B_{1/2}} \L{2}{\nabla \Phi}{B_{1/2}} \to 0.
\end{split}
\]
Moreover if $N>2$ then ${2^*}' =\frac{2N}{N+2}<N$. So by Sobolev embedding we arrive at
\[
\begin{split}
\int_{B_{1/2}} |f_k(x)\Phi(x) |\,dx &\le \|f_k\|_{L^{{2^*}'}(B_{1/2})} \|\Phi \|_{L^{2^*}(B_{1/2})}\\
&\le C(N)\LN{f_k}{B_{1/2}} \L{2^*}{\Phi }{B_{1/2}}\\
&\le  \frac{C(N)}{k}\|\Phi \|_{H^{1}(B_{1/2})}\to 0.\\
\end{split}
\]
In the case where if $N=2$ we have $H^1(B_{1/2})\hookrightarrow L^N(B_{1/2})$, therefore
\[
\begin{split}
\int_{B_{1/2}} |f_k(x)\Phi(x) |\,dx &\le \|f_k\|_{L^{N}(B_{1/2})} \|\Phi \|_{L^{N'}(B_{1/2})}\\
&\le  \frac{C(N)}{k}\|\Phi \|_{H^{1}(B_{1/2})}\to 0.\\
\end{split}
\]
Therefore, the first and third terms in \eqref{e2.} tend to zero as $k\to \infty$. We obtain that for every $\Phi\in H_0^1(B_{1/2})$
\be\label{e0.4}
\lim_{k\to \infty} \int_{B_{1/2}} \big (A(0) \nabla u_k\cdot \nabla \Phi \big ) \,dx = \int_{B_{1/2}}\big (A(0) \nabla u_0\cdot \nabla \Phi \big ) \,dx =0.
\ee
This proves $u_0\in H^1(B_{1/2})$ is a weak solution to $\dive(A(0)\nabla u_0)=0$ in $B_{1/2}$.
\end{proof}
We employ Widman's hole filling technique in the proceeding proposition. The following classical lemma is useful
\begin{lemma}\label{useful lemma} $($c.f. \cite[Lemma 6.1]{eg05}$)$
Let $Z(t)$ be a bounded non negative function defined in $[\rho , r]$. And assume that we have for $\rho \le s <t \le r$
$$
Z(s) \le \theta Z(t) + \frac{A}{|s-t|^2}+C
$$
for some $0\le \theta <1$ and $A,C\ge 0$. Then we have 
$$
Z(\rho) \le C(\theta) \Big [ \frac{A}{|\rho -r|^2 }+C \Big ].
$$
\end{lemma}

\begin{proposition}\label{step 0.15}
Let $u\in H^{1}(B_{1})\cap L^{\infty}(B_1)$ be a weak solution of \eqref{PDE} such that $\Linfty{u}{B_{1}} \le 1$. Then for every $\eps>0$ there exists $0<\delta (\eps, N,\mu) <\frac{1}{2}$ such that if 
$$
\max \Big ( \Linfty{A - A(0)}{B_{1}}, \LN{f}{B_{1}} \Big ) \le \delta
$$
then 
$$
\Linfty{u-h}{B_{1/4}}\le \eps
$$
for some $h\in H^{1}(B_{1/2})$ such that $\Linfty{h}{B_{1/2}}\le 1$ and $h$ is a weak solution to $\dive(A(0)\nabla h)=0$ in $B_{1/2}$.
\end{proposition}

\begin{proof}
We start off by claiming that $\L{2}{\nabla u}{B_{1/2}} \le M(N, \mu)$. 

Indeed, let $s,t>0$ be such that $1/2\le s<t \le 1$ and $\eta \in C_0^{\infty}(B_1)$ such that $0\le \eta \le 1$ and
$$
\eta(x)=
\begin{cases}
1 \;\;x\in B_s\\
0\;\;x\in B_1\setminus B_t.
\end{cases}
$$
we can assume that
\be\label{grad eta}
|\nabla \eta| \le \frac{C(N)}{|s-t|}\;\mbox{ in $B_1$}.
\ee
Now, consider $\vf = \eta u \in H_0^{1}(B_1)$ a test function and from the definition of weak solutions in \eqref{weak solution} we have
\[
\begin{split}
\int_{B_1} \Big ( A(x)\nabla u\cdot \nabla (\eta u) \Big )\,dx = -\int_{B_1} f(x) (\eta u)\,dx.
\end{split}
\]
We expand the integral on the LHS, after rearranging the terms and from ellipticity of $A$ we obtain 
\be\label{f1}
\mu \int_{B_s} |\nabla u|^2\,dx\,\le \,\int_{B_t} \eta \big (A(x)\nabla u\cdot \nabla u\big ) \,dx   = -\int_{B_t} f(x) (\eta u)\,dx- \int_{B_t} u\big (A(x)\nabla u\cdot \nabla \eta  \big )\,dx.
\ee
We observe that since $\nabla \eta = 0$ in $B_s$ 
\be\label{f1.1}
\begin{split}
\int_{B_t} \big | u \big ( A(x) \nabla u\cdot \nabla \eta  \big )\big |\,dx &=\int_{B_t\setminus B_s} \big |u\big ( A(x) \nabla u\cdot \nabla \eta  \big )\big |\,dx\\
&\le \frac{1}{\mu} \L{2}{\nabla u}{B_t\setminus B_s} \L{2}{\nabla \eta }{B_1}.
\end{split}
\ee
Moreover, since $\LN{f}{B_{1}}\le \frac{1}{2}$, $\Linfty{\eta}{B_1}\le 1$ and $\Linfty{u}{B_1}\le 1$ 
\be\label{f1.2}
\begin{split}
\int_{B_1} \Big | f(x) (\eta u)\,dx \Big |& \le  C(N) \LN{f}{B_1} \le C(N).
\end{split}
\ee
using the \eqref{f1.1}, \eqref{f1.2} on the third and fourth term of \eqref{f1}, for some $0<\tau <1$, we arrive to 
\[
\begin{split}
\mu \int_{B_s} |\nabla u|^2\,dx &\le \frac{1}{\mu} \L{2}{\nabla u}{B_t\setminus B_s} \L{2}{\nabla \eta }{B_1} + \Big | \int_{B_t} f(x) (\eta u)\,dx \Big |\\
&\le \frac{C(\mu)}{ 2\tau ^{2}} \L{2}{\nabla u}{B_t\setminus B_s}^{2} + 2\tau ^{2} \L{2}{\nabla \eta }{2} ^2 +C(N).
\end{split}
\]
By taking $\tau = \frac{1}{2}$ and using \eqref{f1} we obtain 
\be\label{f2}
\begin{split}
\int_{B_s} |\nabla u|^2\,dx &\le C_1(\mu)\int_{B_t \setminus B_s} |\nabla u|^2\,dx  + \frac{C_2(N)}{|s-t|^2} +C(N).
\end{split}
\ee
Now, we use the Widman's hole filling procedure. Adding the term $C_1(\mu)\int_{B_s} |\nabla u|^2\,dx$ on both sides of \eqref{f2} and we arrive at 
\be\label{f3}
\begin{split}
\int_{B_s} |\nabla u|^2\,dx &\le \frac{C_1}{C_1+1}\int_{B_t} |\nabla u|^2\,dx + \frac{C_2\cdot (C_1+1)^{-1}}{|s-t|^2}+ \frac{C(N)}{C_1+1}
\end{split}
\ee
Now, from Lemma \ref{useful lemma} with $Z(t):= \int_{B_t}|\nabla u|^2\,dx$ we have
\be\label{f4}
\int_{B_{1/2}} |\nabla u|^2\,dx \le C_3(N,\mu ).
\ee
Once $u$ satisfies \eqref{PDE} in $B_{1}$ we can apply Proposition \ref{step 0.1}.  This way, for a given $\eps>0$, we choose $\delta(\eps,  M, \mu,  N) >0$ in Proposition \ref{step 0.1} which corresponds to $M= C_3(N,\mu)^{1/2}$. Therefore obtain $\delta:= \delta(\eps,\mu,N)$ satisfying the required properties.
\end{proof}
\begin{proposition}[Key lemma for PDE]\label{step 0.2}
Let $u\in H^{1}(B_1)\cap L^{\infty}(B_1)$ be a weak solution to the PDE \eqref{PDE} in $B_1$ with $\Linfty{u}{B_1}\le 1$. Then, for every $0<\alpha<1$, there exists $\delta(N,\mu, \alpha)>0$ and $r_0(N,\mu, \alpha)<1/4$ such that if 
$$
\max\left ( \|A(x) - A(0)\|_{L^{\infty}(B_1)}, \LN{f}{B_1} \right )\le \delta
$$
Then 
$$
\sup_{B_{r_0}}|u-u(0)|\le r_0^{\alpha}.
$$
\end{proposition}

\begin{proof}
Let $\eps>0$ that will be suitably chosen later. By Proposition \ref{step 0.15} there exists $\delta(\eps,N,\mu)>0$ such that
\be\label{s0.2e1}
\|u-h\|_{L^{\infty}(B_{1/4})}<\eps,
\ee 
Where $h$ is a weak solution to $\dive\left  ( A(0)\nabla h \right )=0$ in $B_{1/2}$ and $\Linfty{h}{B_{1/2}}\le 1$. By classical elliptic regularity theory, $h$ satisfies the interior gradient estimates \eqref{int est}. Hence 
\be\label{grad bound}
\Linfty{\nabla h}{B_{1/4}}\le C(N,\mu).
\ee 
Fix $\beta = (1+\alpha)/2<1$, we claim the following
\be\label{s0.2e2}
\sup_{B_r}|h-h(0)| \le C(N,\mu)r^{\beta}\;\;\forall r<1/4.
\ee
Indeed, $h-h(0)$ is a weak solution to the same PDE $\dive \left (A(0)\nabla h \right )=0$ and $\Linfty{h-h(0)}{B_{1/2}}\le 2$ from Proposition \ref{step 0.15}. By mean value theorem and internal gradient estimates \eqref{int est}, if $x\in B_{r}$ and $r<1/4$ we have 
\[
\begin{split}
|h(x)-h(0)| &\le \big |\nabla (h(\theta x)-h(0))\big |\cdot |x|,\,\,(0<\theta<1)\\
&\le 2C(N,\mu)\Linfty{h}{B_{1/2}}r\le C(N,\mu)r^{\beta}
\end{split}
\]
Combining equations (\ref{s0.2e1}) and (\ref{s0.2e2}) we get for $r<1/4$
\be\label{s0.2e3}
\begin{split}
\sup_{B_r}|u(x)-u(0)|&\le \sup_{B_r}\Big ( |u(x)-h(x)|+|h(x)-h(0)|+|h(0)-u(0)| \Big ) \\
&\le 2\eps+C(N,\mu )r^{\beta}.
\end{split}
\ee
We select $r_0(N,\mu,\alpha)<1/4$ such that 
$$
C(N,\mu)r_0^{\beta}\le \frac{r_0^{\alpha}}{3}
$$
that is 
$$
r_0=\Big (  \frac{1}{3C(N,\mu)} \Big )^{1/\beta-\alpha}\le \Big (  \frac{1}{3C(N,\mu)} \Big )^{2/1-\alpha}.
$$
Now, we choose $\eps(N,\mu,\alpha)$ in such a way that 
$$
\eps<\frac{r_0^{\alpha}}{3}.
$$
We see that the choice of $\delta$ depends on $\eps$ and since $\eps$ depends on $N$, $\mu$ and $\alpha$ therefore $\delta$ is actually chosen depending on $N$, $\mu$ and $\alpha$. From \eqref{s0.2e3} we get 
$$
\sup_{B_{r_0}} |u-u(0)| \le r_0^{\alpha}.
$$
\end{proof}

\begin{proposition}\label{step 0.3}
Suppose $u\in H^{1}(B_1)\cap L^{\infty}(B_1)$ is a weak solution to the the PDE \eqref{PDE} in $B_1$ with $\Linfty{u}{B_1}\le 1$. Then, for every $0<\alpha<1$, there exists $\delta(N,\mu, \alpha)>0$ and $C(N,\mu,\alpha)>0$ such that if 
$$
\max\left ( \|A(x) - A(0)\|_{L^{\infty}(B_1)}, \LN{f}{B_1} \right )\le \delta
$$
Then 
$$
\sup_{B_{r}}|u-u(0)|\le C (N,\mu,\alpha)r^{\alpha} \,\,\forall r<r_0,
$$
where $r_0$ is as in Proposition \ref{step 0.2}. Precisely speaking, we have $C(N,\mu,\alpha)= r_0^{-\alpha}$.
\end{proposition} 

\begin{proof}
We argue by scaling and claim that for all $k\in \N$
\be\label{s0.3e1}
 \sup_{B_{r_0^{k}}}|u - u(0)|\le r_0^{k\alpha }.
\ee
It follows readily from Proposition \ref{step 0.2} that (\ref{s0.3e1}) holds for $k=1$. Let us suppose it also holds up to $k_0\in \N$. We prove that (\ref{s0.3e1}) holds also for $k=k_0+1$. In order to do that, we define the following rescaled function
\be\label{rescaled function}
\tilde u(y)=\frac{1}{r_0^{k_0\alpha}}\left ( u(r_0^{k_0}y) - u(0) \right ),\;\;y\in B_1.
\ee
From Remark \ref{rescaling 0.1}, we see that $\tilde u$ satisfies the following PDE
\be\label{rescaled PDE}
\dive( \tilde A(y) \nabla \tilde u ) =\tilde f\;\;\mbox{in $B_1$}
\ee
where $\tilde A$ and $\tilde f$ are given by
\[
\begin{split}
&\tilde A(y):=A(r_0^{k_0}y)\\
&\tilde f(y):=r_0^{k_0(2-\alpha)}f(r_0^{k_0}y).
\end{split}
\]
Now, we verify that the PDE \eqref{rescaled PDE} and its solution $\tilde u$ satisfy the assumptions of Proposition \ref{step 0.2}. Indeed, from \eqref{s0.3e1} for $k=k_0$ and \eqref{rescaled function}, we have 
$$
\sup _{B_{1}}|\tilde u|= r_0^{-k_0 \alpha} \sup_{B_{r_0^{k_0}}}|u-u(0)|\le 1.
$$
Moreover, for $\delta>0$ as in Proposition \ref{step 0.2} we see that
$$
\sup_{B_1}|\tilde A - \tilde A(0)|=  \sup_{B_{r_0^{k_0}}} |A-A(0)| \le \delta.
$$
Additionally
$$
\|\tilde f\|_{L^{N}(B_1)}=r_0^{k_0(1-\alpha)}\|f\|_{L^N(B_{r_0^{k_0}})}\le  \delta.
$$
So we can apply the Proposition \ref{step 0.2} for $\tilde A$, $\tilde f$ and $\tilde u$. Hence we obtain 
$$
\sup_{B_{r_0}}|\tilde u|\le r_0^{\alpha}.
$$
Rescaling back to $u$ we arrive at 
$$
\sup_{B_{r_0^{k_0+1}}}|u-u(0)|\le r_0^{(k_0+1)\alpha}.
$$
So the claim (\ref{s0.3e1}) is proven for all $k\in \N$. To finish the proof of Proposition \ref{step 0.3}, let $r\in (0,r_0)$ be given then we can find $k\in \N$ such that $r_0^{k+1}\le r<r_0^k$. From (\ref{s0.3e1}), we see that 
\be\label{explicit}
\sup_{B_r}|u-u(0)|\le \sup_{B_{r_0^{k}}}|u-u(0)|\le r_0^{k\alpha}=r_0^{(k+1)\alpha}\frac{1}{r_0^{\alpha}}\le {r_0^{-\alpha}}r^{\alpha}.
\ee
Finally we can take $C(N,\mu,\alpha) = \frac{1}{r_0^{\alpha}}$. Then  Proposition \ref{step 0.3} is proven.
\end{proof}

In the sequel, we remove the small regime condition on the previous results. As a matter of fact we prove the main result of this section.

\begin{proposition}\label{step 0.4}
Suppose $u\in H^{1}(B_1)\cap L^{\infty}(B_1)$ be a weak solution to \eqref{PDE} in $B_1$. Then for any $0<\alpha<1$, we have $u\in C^{\alpha}(B_{1/2})$ with the following estimates
\be\label{est}
[u]_{C^{\alpha}(B_{1/2})} \le   C(N,\alpha, \mu,\o_{A,\overline{B_{3/4}}})  \left ( \Linfty{u}{B_{1}}+\LN{f}{B_{1}} \right ).
\ee
Here $\o_{A,\overline{B_{3/4}}}$ is uniform the modulus of continuity of $A$ in $\overline{B_{3/4}}$ (c.f. Definition \ref{mod of cont}). In particular, we have 
\be\label{est2}
\|u \|_{C^{\alpha}(B_{1/2})} \le C(N,\alpha,\mu,\o_{A,\overline{B_{3/4}}})  \left ( \Linfty{u}{B_{1}}+\LN{f}{B_{1}}\right ).
\ee
\end{proposition}

\begin{proof}
Since $A\in C^{N\times N}(B_1)$, then the coefficient matrix $A$ is uniformly continuous in $\overline{B_{3/4}}$. Hence, there exists $s_0:= s_0(\o_{A,\overline{B_{3/4}}}, N, \mu, \alpha)$ such that $0<s_0<1/4$ and for every $x_0\in B_{1/2}$ we have 
\be\label{smallness}
\sup_{B_{s_0}(x_0)} |A-A(x_0)|\le \delta(N,\mu,\alpha)
\ee
for $\delta(N,\mu,\alpha)>0$ as in Porposition \ref{step 0.3}. We fix now $x_0\in B_{1/2}$ and define the following rescaled function 
\be\label{w}
w(y):= \frac{u(x_0 + s_0 y)}{\left ( \Linfty{u}{B_1}+\frac{s_0}{\delta}\LN{f}{B_1} \right )},\,\,\;y\in B_1.
\ee
By Remark \ref{rescaling 0.1}, we can easily verify that the function $w\in W^{1,p}(B_1)$ satisfies the following PDE
\be\label{PDE2}
\dive(\b{A}(y)\nabla v)  = \b{f}\,\,\mbox{in $B_1$}.
\ee
Where $\b{A}$ and $\b{f}$ are defined as follows
\[
\begin{split}
\b{A}(y)&:= A(x_0+s_0 y)\\
\b{f}(y)&: = \frac{s_0^2\cdot f(x_0+s_0y)}{\left ( \Linfty{u}{B_1}+\frac{s_0}{\delta}\LN{f}{B_1} \right )},\;y\in B_1.
\end{split}
\]
From \eqref{w}, we can see that 
\be\label{Linfty w <=1}
\Linfty{w}{B_1}\le 1.
\ee
We rewrite \eqref{smallness} as follows 
\be\label{smallness2}
\sup_{B_{s_0}(x_0)} |A-A(x_0)|=\Linfty{\b{A}-\b{A}(0)}{B_1}\le\delta(N,\mu,\alpha).
\ee 
Moreover, from the definition of $\b{f}$, we see that 
\be\label{smallness3}
\begin{split}
\LN{\b{f}}{B_1}&= \frac{s_0\cdot \LN{f}{B_{s_0}(x_0)}}{\left ( \Linfty{u}{B_1}+\frac{s_0}{\delta}\LN{f}{B_1)} \right )} \le \delta(N,\mu,\alpha).
\end{split}
\ee
From \eqref{Linfty w <=1}, \eqref{PDE2}, \eqref{smallness2} and \eqref{smallness3} we see that $w, \b{A}$ and $\b{f}$ satisfy the assumptions of Proposition \ref{step 0.3}, and therefore from Proposition \ref{step 0.3} we have 
$$
\sup_{B_r}|w-w(0)|\le \frac{1}{r_0^{\alpha}} r^{\alpha},\;\;\forall r\le r_0(N,\mu,\alpha)
$$
Rescaling back to $u$, we obtain for every $x_0\in B_{1/2}$
$$
\sup_{z\in B_{s_0r}(x_0)} \Bigg | \frac{u-u(x_0)}{\left ( \Linfty{u}{B_1}+\frac{s_0}{\delta}\LN{f}{B_1} \right )} \Bigg | \le C r^{\alpha},\;\;r\le r_0.
$$
Therefore
\be\label{est1}
\sup_{B_r(x_0)}|u-u(x_0)| \le \frac{1}{r_0^{\alpha}s_0^{\alpha}} \left ( \Linfty{u}{B_1}+\frac{s_0}{\delta}\LN{f}{B_1} \right )r^{\alpha},\;\;\forall\, r\le s_0 r_0.
\ee
Let $x,y\in B_{1/2}$. If $|x-y|\le s_0r_0$, then we replace $x$ by $x_0$ and take $r=|x-y|$ in the above equation and we have 
$$
|u(x)-u(y)|\le \frac{1}{r_0^{\alpha}s_0^{\alpha}} \left ( \Linfty{u}{B_{1}}+\frac{s_0}{\delta}\LN{f}{B_{1}} \right ) |x-y|^{\alpha}.
$$
Otherwise, if $|x-y|> s_0r_0$ then $\frac{|u(x)-u(y)|}{|x-y|^{\alpha}} \le \frac{2\Linfty{u}{B_{1}}}{s_0^{\alpha}r_0^{\alpha}} $. In other words,
$$
|u(x)-u(y)|\le \frac{2\Linfty{u}{B_{1}}}{s_0^{\alpha}r_0^{\alpha}}|x-y|^{\alpha}.
$$
Finally combining the two estimates above we have for all $x,y\in B_{1/2}$
\be\label{refeq}
\begin{split}
|u(x)-u(y)| &\le\frac{C}{r_0^{\alpha}s_0^{\alpha}}   \left ( \Linfty{u}{B_{1}}+\frac{s_0}{\delta}\LN{f}{B_{1}} \right ) |x-y|^{\alpha}\\
&\le \frac{C(1+\frac{1}{\delta})}{r_0^{\alpha}s_0^{\alpha}}\left ( \Linfty{u}{B_1}+\LN{f}{B_1} \right )|x-y|^\alpha.
\end{split}
\ee
This proves Proposition \ref{step 0.4}.
\end{proof}

\begin{remark}\label{no modulus rem}
The dependence of constant at right hand side of \eqref{est} and \eqref{est2} on the modulus of continuity $\o_{A, \overline{B_{3/4}}}$ can be removed as long as we are under the smallness condition on the oscillations of the coefficient $A$. That is if we have 
\be\label{refeq1}
\Linfty{A-A(0)}{B_1}\le \delta
\ee
for $\delta (N,\mu,\alpha)$ as in Lemma \ref{step 0.3}. Then under this condition we have 
\be\label{no modulus eq}
\|u \|_{C^{\alpha}(B_{1/2})} \le C(N,\mu,\alpha)  \left ( \Linfty{u}{B_{1}}+\LN{f}{B_{1}} \right ).
\ee
To see this, we observe that in \eqref{refeq}, the dependence on $\o_{A,\overline{B_{3/4}}}$ is coming from $s_0$. This positive number $s_0$ is chosen so small that $\Linfty{A-A(0)}{B_{s_0(x_0)}}\le \delta$. This step is not required if we already assume \eqref{refeq1}.

\end{remark}

Proposition \ref{step 0.4} is scaled version of the regularity estimates for weak solution to the PDE \eqref{PDE}. From the classical rescaling argument, we can also prove the following result for any general domain $\O$.

\begin{corollary}\label{step 0.45}
Let $u\in H^{1}(B_\rho(x_0))\cap L^{\infty}(B_{\rho}(x_0))$ be a weak solution to
$$
\dive\left (A(x)\nabla u\right )=f\;\;\mbox{in $B_\rho(x_0)$}
$$
where $A\in C(B_\rho(x_0))^{N\times N}$ and there exists $0<\mu<1$ such that $\mu |\xi|^2 \le \< A(x)\xi ,\xi \>\le \frac{1}{\mu}|\xi|^2$ a.e. in $B_\rho(x_0)$ and for all $\xi \in \R^N$. $f\in L^N(B_\rho(x_0))$. Then for any $0<\alpha<1$, $u\in C^{0,\alpha}(B_{\frac{\rho}{2}}(x_0))$ and if $0<\rho<1$ we have
\be\label{estg1}
\|u\|_{C^{0,\alpha}(B_{\frac{\rho}{2}}(x_0))} \le   \frac{C(N,\mu,\alpha,\o_{A,\overline{B_{\frac{3\rho}{4}}(x_0)}})}{\rho^{\alpha}}  \left ( \Linfty{u}{B_{\rho}(x_0)}+\rho \LN{f}{B_{\rho}(x_0)} \right ).
\ee
Here $\o_{A,\overline{B_{\frac{3\rho}{4}}(x_0)}}$ is the uniform modulus of continuity of the coefficient $A$ in $\overline{B_{\frac{3\rho}{4}}(x_0)}$. Furthermore, the dependence on $\o_{A,\overline{B_{\frac{3\rho}{4}}(x_0)}}$ can be dropped under the smallness assumption  
\be\label{refeq2}
\Linfty{A-A(x_0)}{B_\rho(x_0)}\le \delta 
\ee
for $(\delta$ as in Lemma \ref{step 0.3} $)$. Then if $0<\rho<1$
\be\label{estg2}
\|u\|_{C^{0,\alpha}(B_{\frac{\rho}{2}}(x_0))} \le   \frac{C(N,\mu,\alpha)}{\rho^{\alpha}}   \left ( \Linfty{u}{B_{\rho}(x_0)}+\rho \LN{f}{B_{\rho}(x_0)} \right ).
\ee
\end{corollary}
\begin{proof}
The proof of Corollary \ref{step 0.45} follows readily by scaling. In fact, we just consider the function defined as $v(y):= \frac{u(x_0+\rho y)}{\rho}$ and apply Proposition \ref{step 0.4} together with Remark \ref{no modulus rem}. Indeed, from Remark \ref{rescaling 0.1}, $v\in H^1(B_1)\cap L^{\infty}(B_1)$ is a weak solution to the PDE $\dive(\bar A(y)\nabla v) = \bar f$ in $B_1$. Where $\bar A$ and $\bar f$ are given by 
$$
\bar A(y) = A(x_0+\rho y),\;\;\bar f(y) = \rho f(x_0+\rho y).
$$
We can obtain the estimates from Proposition \ref{step 0.4} on $v$ and by translating information back to $u$.
\end{proof}
We prove the main result of this section
\section{Proof of Theorem \ref{PDE regularity}}

\begin{proof}[Proof of Theorem \ref{PDE regularity}]
Let $D\subset \subset B_1$. We observe that $D\subset \subset B_{d} \subset \subset B_1$ where $d=1-\frac{\dist(D,\partial B_1)}{2}$. Since $A\in C(B_1)^{N\times N}$, then the matrix $A$ is uniformly continuous in $\overline{B_d}$. Thus, there exists a modulus of continuity $\o_{A, \overline {B_d}}$ (c.f. Definition \ref{mod of cont}) given by
$$
 \o_{A,\overline{B_d}}(t) := \begin{cases} \Bigg ( \sup_{\substack{ |x-y|<t\\x,y\in \overline{B_d}}}|A(x)-A(y)|\Bigg ) \qquad t \le 2d\\
 \o_{A,\overline{B_d}}(2d)\qquad\qquad\qquad\qquad \;\;\;\;\;t>2d .
 \end{cases}
$$
{We set $t_0$ as  
$$
t_0 := t_0(\o_{A,\overline{B_d}},\delta) = \sup \Big \{ t \;\Big | \;\o_{A,\overline{B_d}}(t)\le \delta \Big \}
$$ 
as well as $s_0 := \min\Big (t_0,\frac{\dist(D,\partial B_1)}{4} \Big )$. Since $\o_{A,\overline {B_d}}$ is a non-decreasing function we have $\o_{A,\overline{B_d}}(s_0)\le \delta$. }Now since 
$$
D\subset \bigcup_{x\in D}B_{s_0}(x) \subset \overline {B_d}
$$
we have 
\be\label{smallness1..}
\sup_{B_{s_0}(x)}|A-A(x)| \le \o_{A,\overline{B_d}}(s_0) \le {\delta},\;\forall x\in D.
\ee
We note that $u$ is also an absolute minimizer of $J(\cdot;B_{s_0}(x))$, by \eqref{smallness1..} and \eqref{estg2} in Corollary \ref{step 0.45} we have for all $y\in B_{s_0/2}(x)\cap D$ 
\be\label{near x_0..}
|u(x)-u(y)|\le \frac{C(N,\alpha, \mu)}{s_0^{\alpha}}\left (\Linfty{u}{B_{1}}+s_0\LN{f}{B_{1}}  \right )|x-y|^{\alpha}.
\ee
Now, if $x, y\in D$ are such that $|x-y| \ge s_0/2$, then 
\be\label{far from x_0..}
{|u(x)-u(y)|}\le 2^{1+\alpha} \frac{\Linfty{u}{B_1}}{s_0^{\alpha}}|x-y|^{\alpha}.
\ee
By combining \eqref{near x_0..} and \eqref{far from x_0..}, we arrive to $\Big ($since $s_0\le\frac{\dist(D,\partial B_1)}{4}\le\frac{\diam(B_1)}{4} =\frac{1}{2}<1$$\Big )$
\be\label{precise est..}
[u]_{C^{\alpha}(D)} \le \frac{C(N,\alpha, \mu)}{s_0^{\alpha}}\left ( \Linfty{u}{B_{1}}+\LN{f}{B_{1}} \right ).
\ee
{We observe from the definition of $s_0$ that 
\be\label{precise est2..}
[u]_{C^{\alpha}(D)} \le \begin{cases}
\frac{C(N,\alpha, \mu)}{t_0^{\alpha}}\left ( \Linfty{u}{B_{1}}+\LN{f}{B_{1}} \right ) \;\;\;\;\;\;\;\;\;\;\;\;\;\;\;\;\;\;\;\mbox{if $t_0\le \frac{\dist(D,\partial B_1)}{4}$}\\[7pt]
\frac{4^\alpha\cdot C(N,\alpha, \mu)}{\dist(D,\partial B_1)^\alpha}\left ( \Linfty{u}{B_{1}}+\LN{f}{B_{1}} \right )\;\; \;\;\;\;\;\;\;\;\;\;\;\;\;\mbox{if $t_0\ge \frac{\dist(D,\partial B_1)}{4}$.}
\end{cases}
\ee
}
{In order to control the first term in the equation above by a universal multiple of $\dist(D,\partial B_1)^{-\alpha}$, we observe that} once $t_0>0$ depends only on the modulus of continuity $\o_{A, \overline{B_d}}$ and $\delta$, there exists a $n_0:= n_0(\o_{A, \overline{B_d}},\delta)=n_0(N,\mu,\alpha,\o_{A,\overline{B_d}})\in \N$ depending only on $t_0$ such that  $\frac{2}{n_0}\le t_0$. Hence
\be\label{s0}
\frac{\dist(D,\partial B_1)}{n_0}\le \frac{\diam(B_1)}{n_0} =\frac{2}{n_0}\le t_0 \implies \frac{1}{t_0^{\alpha}}\le \frac{n_0^{\alpha} }{\dist(D,\partial B_1)^{\alpha}}= \frac{C(N,\mu,\alpha,\o_{A,\overline{B_d}}) }{\dist(D,\partial B_1)^{\alpha}}.
\ee
\begin{comment}
Thus the estimate \eqref{precise est} becomes
\be\label{last est}
[u]_{C^{\alpha}(D)} \le \frac{C(N, \alpha, \mu, q_2\lambda_+, \o_{A_{\pm}, D})} {\dist(D,\partial B_1)^{\alpha}}\left ( 1+\Linfty{u}{B_1}+\LN{f}{B_1} \right ).
\ee
We consider a finite cover of $\overline D\subset \subset B_1$ 
$$
\overline D \subset \bigcup_{i=1}^n \Big \{ B_{\frac{s_0}{8}}(x_i): \,x_i\in D \Big \}.
$$
Let $x,y\in D$, and assume $x\in B_{\frac{s_0}{8}(x_i)}$ for some $i\in \{1,...,n\}$. If $|x-y|\le \frac{s_0}{16}$, then 
$$
|y-x_i|\le |x-y|+|x-x_i|\le \frac{s_0}{16}+\frac{s_0}{8}= \frac{3s_0}{19}<<\frac{s_0}{2}.
$$
Therefore, $x,y\in B_{s_0/2}(x_i)$. From \eqref{last est} with $x_0=x_i$ we obtain
\be\label{x,y near}
|u(x)-u(y)| \le \frac{C(N, \alpha, \mu, q_2\lambda_+, \o_{A_{\pm}, D})} {\dist(D,\partial B_1)^{\alpha}}\left ( 1+\Linfty{u}{B_1}+\LN{f}{B_1} \right )|x-y|^{\alpha}
\ee
If $|x-y| \ge \frac{s_0}{8}$ then from \eqref{s0} we have
\be\label{x,y far}
|u(x)-u(y)| \le \frac{8^{\alpha}\cdot 2\Linfty{u}{B_1}}{s_0^{\alpha}}|x-y|^{\alpha}\le \frac{C(\o_{A_{\pm},D},\alpha)\Linfty{u}{B_1}}{\dist(D,\partial B_1)^{\alpha}}|x-y|^{\alpha}.
\ee
\end{comment}
Now, \eqref{precise est2..} becomes
\be\label{general est..}
[u]_{C^{\alpha}(D)} \le \frac{C(N,\alpha, \mu, \o_{A_{\pm},\overline{B_d}})}{\dist(D,\partial B_1)^{\alpha}}\left ( \Linfty{u}{B_{1}}+\LN{f}{B_{1}} \right ).
\ee
To finish we observe that Theorem \ref{PDE regularity} follows from \eqref{general est..} by taking $D$ as $B_r$ for any $r<1$, once $\dist(D,\partial B_1)^{\alpha} = \dist(B_r,\partial B_1)^{\alpha} = (1-r)^{\alpha}$ and $d=\frac{1+r}{2}$. 
\end{proof}

\section{Approximation lemma via compactness}\label{compactness}
Our strategy is to make use of regularity theory of limiting functional and prove regularity of minimizers. We were inspired from the method used in \cite[Lemma 5.4]{AC81}. First we show that if the functions $A_+$ and $A_-$ are very close to different constants, the regularity of a minimizer is close to $C^{0,1^-}$.
\begin{remark}\label{lim of seq}
Although elementary, we make a systematic use of limsup and liminf properties (especially $*$ and $**$) in our proof of compactness. For the reader's convenience, we record them here. For any  bounded sequences of real numbers $a_k$ and $b_k$, the following properties hold (for the proofs check for instance \cite[Theorem 3.127]{D17})
\begin{enumerate}[label=\textbf{($\II$\arabic*)}]
\item \label{I1} \label{I2}  $ \limsup_{k\to \infty}(-a_k)= -\liminf_{k\to \infty} (a_k) $ and $ \liminf_{k\to \infty}(-a_k)= -\limsup_{k\to \infty} (a_k) $.
\item \label{I3} \label{I4} Also we have 
\[
\begin{split}
\liminf_{k\to \infty}a_k + \liminf_{k\to \infty}b_k &\stackrel{*}{\le} \liminf_{k\to \infty}(a_k+b_k)\\
&\le \limsup_{k\to \infty}a_k + \liminf_{k\to \infty}b_k\\
&\le \limsup_{k\to \infty}(a_k+b_k)\\
&\stackrel{**}{\le} \limsup_{k\to \infty}a_k+\limsup_{k\to \infty}b_k.
\end{split}
\]
\end{enumerate}
\end{remark}

\begin{remark}
In the forthcoming proofs, we consider the functional $J(\cdot;B_1):=J_{A,f,Q}(\cdot; B_1)$,  
where $B_1:= B_1(0)$ and
$$
J(\cdot;B_1):=J_{A,f,Q} (v;B_1) =\int_{B_1} A(x,u)|\nabla u|^2+f(x,u)u+ Q(x)\lambda (v)\,dx.
$$
Up to rescaling, $B_1$ represents a small ball contained inside a general domain $\O$.  We already know that local minimizers of $J(\cdot;\O)$ are locally bounded in $\O$ c.f. \cite[Theorem 7.3]{eg05},  \cite[Theorem 2.3]{HS21}. Therefore it is reasonable to assume that for a minimizer $u\in H^1(B_1)$ of $J_{A,f,Q} (u;B_1)<\infty $ and $\|u\|_{L^{\infty}(B_1)} \le 1$.  These assumptions do not compromise any generality.
\end{remark}

We now approximate minimizers by regular ones under the small regime scenario.{
\begin{theorem} \label{step 1}
Let $u\in H^{1}(B_{1/2})$ be an absolute minimizer of $J(\cdot, B_{1/2})$ such that $J(\cdot; B_{1/2})$ satisfies structural conditions \ref{G1}-\ref{G5}. Assume also that $\L{2}{\nabla u}{B_{1/2}}\le M$ and

\noindent $\Linfty{u}{B_{1/2}}\le 1$. Then for any $\eps>0$ there exists $\delta(\eps ,M, \mu, q_2,\lambda_+, N)>0$ such that if 
$$
 \max \Big (\|A_{\pm}-A_{\pm}(0)\|_{L^{\infty}(B_{1/2})},\, \|f_{\pm}\|_{L^N(B_{1/2})},\, \|Q-Q(0)\|_{L^{\infty}(B_{1/2})} \Big ) \le \delta 
$$
then
$$
\|u-h\|_{L^{\infty}(B_{1/4})}\le\eps
$$
where $h$ is an absolute minimizer of $$J_0(v; B_{1/2}):= \int_{B_{1/2}}\Big ( A(0,v)|\nabla v|^2 + Q(0)\lambda(v) \Big )\,dx.$$
Moreover, $\Linfty{h}{B_{1/2}}\le 1$.
\end{theorem}
}
\begin{proof}
{
Let us suppose for the sake of contradiction that the statement of Theorem \ref{step 1} is not true. This implies that there exists an $\eps_0>0$ and sequences $A_k, f_k, Q_k$ such that $A_{\pm,k}\in C(B_{1/2})$, $Q_k\in L^{\infty}(B_{1/2})$ and $f_k\in L^{N}(B_{1/2})$ with $\|A_{\pm,k}-A_{\pm}(0)\|_{L^{\infty}(B_{1/2})}<\frac{1}{k}$, $\|f_{\pm,k}\|_{L^N(B_{1/2})}<\frac{1}{k}$ and $\|Q_k-Q(0)\|_{L^{\infty}(B_{1/2})} < \frac{1}{k}$ as well as corresponding absolute minimizers $u_k\in H^1(B_{1/2})$ of the functional $J_k(\cdot;B_{1/2})$ (defined in \eqref{def of Jk} below) such that $\L{p}{\nabla u_k}{B_{1/2}}\le M$ and $\| u_k \|_{L^{\infty}(B_{1/2})} \le 1$ also satisfy
\be\label{absurdum}
\|u_k-h\|_{L^{\infty}(B_{1/4})}>\eps_0
\ee
for all $h\in H^1(B_{1/2})\cap L^{\infty}(B_{1/2})$ with $\Linfty{h}{B_{1/2}}$ that are absolute minimizers of the functional of the type 
$$
J_0(v;B_{1/2}):= \int_{B_{1/2}}\Big ( A(0,v) |\nabla v|^2 +Q(0)\lambda(v) \Big )\,dx.
$$
The functional $J_k(\cdot;B_{1/2})$ is given by
\be\label{def of Jk}
J_k(v;B_{1/2}):= \int_{B_{1/2}} \Big (A_k(x,v)|\nabla v|^2+f_k(x,v)v +Q_k(x)\lambda (v) \Big )\,dx.
\ee
We recall that $\lambda(s)=\lambda_+\chi_{\{s>0\}}+\lambda_-\chi_{\{s\le 0\}}$ and
\be\label{def of coeff}
A_k(x,v) = A_{+,k}(x)\chi_{\{v>0\}}+A_{-,k}(x)\chi_{\{v\le0\}},\; f_k(x,v)= f_{+,k}(x)\chi_{\{v>0\}}+f_{-,k}(x)\chi_{\{v\le0\}}.
\ee
From \ref{G1}-\ref{G5}, it follows that the functional $J_k(\cdot;B_{1/2})$ and its integrand satisfy the structural condition in \cite[equation 7.2]{eg05} uniformly in $k$. Therefore, from \cite[Theorem 7.6]{eg05} it follows the existence of $\beta_0:= \beta_0 ( \mu, q_2,\lambda_+,N)\in (0,1)$ and $C_0:= C_0 (\mu, q_2,\lambda_+,N)>0$ such that for every $k\in \N$ we have
$$
\| u_k\|_{C^{\beta_0}(\overline{B_{1/4}})}  \le C_0.
$$
By Arzela Ascoli theorem, there exists $u_0\in C^{0,\beta_0}(B_{1/4})$ such that 
\be\label{uniform convergence p=2} 
u_k\to u_0 \mbox{ in $L^{\infty}(B_{1/4})$ up to a subsequence. }
\ee
Since $\L{2}{\nabla u_k}{B_{1/2}}\le M$ and $\L{2}{u_k}{B_{1/2}} \le C(N)\cdot\Linfty{u_k}{B_{1/2}}\le C(N)$ for every $k$, $u_k$ is a bounded sequence in $H^1(B_{1/2})$.  That is 
\be\label{bound on uk}
\| u_k \|_{H^1(B_{1/2})} \le C(M, N).
\ee
Therefore, $u_k$ converges weakly to $u_0$ (up to a subsequence)
$$
u_k \wto u_0 \;\;\mbox{weakly in $H^1(B_{1/2})$}.
$$
From \eqref{bound on uk} and weak lower semicontinuity of $H^1(B_{1/2})$ norm, this implies 
\be\label{bound on u0}
\|u_0\|_{H^1(B_1)} \le C(M,N).
\ee
Also $\Linfty{u_0}{B_{1/2}}\le 1$ from Lemma \ref{weak convergence preserves inequality}. 
In order to simplify the steps to come, we introduce $\GG_k(\cdot;B_{1/2})$ as follows
$$
\GG_k(v;B_{1/2}) := \int_{B_{1/2}} \Big ( (A_k(x,v) -A(0,v) )|\nabla v|^2 +f_k(x,v)v +(Q_k(x)-Q(0)) \lambda(v) \Big )\,dx.
$$
We observe that for any function $v\in H^1(B_{1/2})$, we have the following splitting
\be\label{add subtract}
J_k(v;B_{1/2}) = J_0(v;B_{1/2}) +\GG_k(v;B_{1/2}).
\ee
Now we make a claim
\be\label{M}\tag{$\pmb{\MM}$}
\mbox{\textbf{Claim:} $u_0$ is an absolute minimizer of $J_0(\cdot;B_{1/2})$.}
\ee
Proving the claim \eqref{M} will lead us to a contradiction to \eqref{absurdum} and thus to the proof of Theorem \ref{step 1}.

In order to prove the claim \eqref{M}, let us take a function $v\in H^1(B_{1/2})$ be such that $v-u_0\in H_0^1(B_{1/2})$. To ease the notation, we relabel the constants $A_{\pm}(0)$ as follows
\[
\begin{split}
A_+(0)= : a^2\\
A_-(0)= : b^2.
\end{split}
\]
We observe that $\mu \le a^2, b^2 \le \frac{1}{\mu}$ (c.f. \ref{G3}).  For $r\in (0,\frac{1}{2})$ we set $\eta_{r}\in C_c^{\infty}(\R^N)$ so that
\be\label{eta}
\eta_r=1 \;\;\mbox{in $B_{\frac{1}{2}-r}$, \qquad  $\supp(\eta_r)\subset \overline {B_{1/2}}$,  \qquad$|\nabla \eta| \le \frac{C(N)}{r}$.}
\ee
Now, we consider $v_{k,r} \in H^1(B_{1/2})$ such that 
\be\label{def of vk}
\TT_{a,b} (v_{k,r}) := \TT_{a,b}(v) + (1-\eta_{r}) (\TT_{a,b}(u_k) - \TT_{a,b}(u_0)).
\ee
The function $v_{k,r}\in H^1(B_{1/2})$ above is well defined due to the fact  that the $\TT_{a,b}$-operator is a bijective map in $H^{1}(B_{1/2})$ (c.f. Lemma \ref{one one onto}-\ref{q1}).  

\vspace{0.3cm}
Let us verify that $v_{k,r}$ is an admissible competitor for the minimality of $u_k$ for the functional $J_k(\cdot;B_{1/2})$. Indeed, we observe that since $\eta_r\in C_c^{\infty}(B_{1/2}) \subset H_0^{1}(B_{1/2})$ then 
$$
\Tr \left ( \eta_r \left ( \TT_{a,b}(u_k) - \TT_{a,b}(u_0) \right ) \right ) =0\mbox{ on $\partial B_{1/2}$.}
$$ 
Now, we have
\vspace{0.3cm}
\[
\begin{split}
\Tr\left (\TT_{a,b}(v_{k,r}) \right ) &= \Tr(\TT_{a,b}(v)) + \Tr \left ( (1-\eta_r)(\TT_{a,b}(u_k) - \TT_{a,b}(u_0)) \right )\\
&= \Tr \left ( \TT_{a,b}(v) \right ) + \Tr \left ( \TT_{a,b}(u_k) - \TT_{a,b}(u_0) \right ) - \Tr \left ( \eta_r \left ( \TT_{a,b}(u_k) - \TT_{a,b}(u_0) \right ) \right )\\
&= \Tr \left ( \TT_{a,b}(v) \right ) + \Tr \left ( \TT_{a,b}(u_k) - \TT_{a,b}(u_0) \right ) = \Tr \left ( \TT_{a,b}(v) \right )+ \Tr \left ( \TT_{a,b}(u_k) \right ) - \Tr \left ( \TT_{a,b}(u_0) \right ).\\ 
\end{split}
\]
Using Proposition \ref{prop of Tab}-\ref{p5}, and the fact that $\Tr(u_0) = \Tr(v)$ on $\partial B_{1/2}$ we arrive at
\[
\begin{split}
\TT^{\partial B_{1/2}}_{a,b}\big (\Tr (v_{k,r}) \big ) = \TT^{\partial B_{1/2}}_{a,b}\big (\Tr (v) \big ) + \TT^{\partial B_{1/2}}_{a,b}\big (\Tr (u_k) \big )- \TT^{\partial B_{1/2}}_{a,b}\big (\Tr (u_0) \big ) = \TT^{\partial B_{1/2}}_{a,b}\big (\Tr (u_k) \big ).
\end{split}
\]
Since $\TT_{a,b}^{\partial B_{1/2}}$ is a bijective map (c.f.  Proposition \ref{one one onto}-\ref{q1.5}), we conclude
\be\label{admissible competitor}
\Tr(v_{k,r}) = \Tr(u_k)\mbox{ on $\partial B_{1/2}$.}
\ee
Therefore, by minimality of $u_k$ for $J_k(\cdot;B_{1/2})$
\be\label{minimality of uk}
J_k(u_k;B_{1/2}) \le J_k(v_{k,r};B_{1/2}).
\ee
From the splitting in \eqref{add subtract}, we can write
\be\label{min of uk 1}
J_0(u_k;B_{1/2}) + \GG_{k}(u_k;B_{1/2}) \le J_0(v_{k,r};B_{1/2})+\GG_k(v_{k,r};B_{1/2}).
\ee
In order to further ease the notation, we relabel $\TT_{a,b}(v_{k,r}) $, $\TT_{a,b}(v) $, $\TT_{a,b}(u_0) $, $\TT_{a,b}(u_k) $ as follows
\be\label{defs}
\begin{split}
\TT_{a,b}(v_{k,r}) = :V_{k,r},\;\;
\TT_{a,b}(v) = : V,\;\;
\TT_{a,b}(u_0)= :U_0,\;\;
\TT_{a,b}(u_k)= : U_k .
\end{split}
\ee
We can rewrite the definition of $v_{k,r}$ in \eqref{def of vk} as
\be\label{def of  Vkr}
V_{k,r} = V + (1-\eta_r )(U_k -U_0).
\ee
Next, we show that $\|v_{k,r}\|_{H^1(B_{1/2})}$ and $\|V_{k,r}\|_{H^1(B_{1/2})}$ are bounded sequences independent of $k$. First we observe that from Proposition \ref{one one onto}-\ref{q2}, \eqref{bound on uk} and \eqref{bound on u0} we have
\be\label{bounds on Uk and U0}
\begin{split}
\|U_k\|_{H^1(B_{1/2})} &\le \max(a,b) \cdot \|u_k\|_{H^1(B_{1/2})} \le \frac{1}{\sqrt{\mu}}C(N,M)\\
\|U_0\|_{H^1(B_{1/2})} &\le \max(a,b) \cdot\|u_0\|_{H^1(B_{1/2})} \le \frac{1}{\sqrt{\mu}}C(N,M)\\
\|V\|_{H^1(B_{1/2})} &\le \max(a,b)\cdot \|v\|_{H^1(B_{1/2})} \le \frac{1}{\sqrt{\mu}}\|v\|_{H^1(B_{1/2})}.
\end{split}
\ee
In order to bound $\int_{B_{1/2}}|V_{k,r}|^2\,dx$ we observe that from \eqref{bounds on Uk and U0} we have 
\[
\begin{split}
\int_{B_{1/2}} |V_{k,r}|^2\,dx &= \int_{B_{1/2}}\big |V+(1-\eta_r)(U_k-U_0) \big |^2\,dx\\
&\le 2 \int_{B_{1/2}}\big ( |V|^2+|U_k|^2+|U_0|^2 \big )\,dx\\
&\le \frac{1}{\mu} \Big ( \|v\|^2_{H^1(B_{1/2})}  + C(N,M)\Big ) = C(N,M,\mu) \Big ( \|v\|^2_{H^1(B_{1/2})}+1 \Big ).
\end{split}
\]
Now, we estimate $\int_{B_{1/2}}|\nabla V_{k,r}|^2\,dx$. Again using \eqref{bounds on Uk and U0} and \eqref{eta} we have 
\be\label{bound on Vkr}
\begin{split}
\int_{B_{1/2}}|\nabla V_{k,r}|^2\,dx&= \int_{B_{1/2}} \big | \nabla \big ( V+(1-\eta_r)(U_k-U_0) \big ) \big |^2\,dx\\
&\le 2\Bigg [ \int_{B_{1/2}}|\nabla V|^2\,dx + \frac{C(N)}{r^2} \int_{B_{1/2}} |U_0 - U_k|^2 \,dx + \int_{B_{1/2}} | \nabla (U_0 -U_k)|^2\,dx\Bigg ]\\
&\le 4 \Bigg [ \int_{B_{1/2}} \Big (  |\nabla V|^2 + |\nabla U_0|^2 + |\nabla U_k|^2 \Big )\,dx  +\frac{C(N)}{r^2} \int_{B_{1/2}} \Big (  |U_0|^2 + |U_k|^2\Big )\,dx\Bigg ]\\
&\le 4 \Bigg [  \frac{1}{\mu} \|v\|^2_{H^1(B_{1/2})} + \frac{1}{\mu}C(N,M) + \frac{1}{r^2 \mu}C(N,M) \Bigg]\\
&\le \frac{\bar C(N,M)}{r^2 \mu} \Big ( \|v\|^2_{H^1(B_{1/2})}+1 \Big )= C_1(N,M,\mu,r)\Big ( \|v\|^2_{H^1(B_{1/2})}+1 \Big ).
\end{split}
\ee
Thus, $r>0$ we have a bound on $\|V_{k,r}\|_{H^1(B_{1/2})}$ independent of $k$. Again from Proposition \ref{one one onto}-\ref{q2} we have 
\be\label{bound on vkr}
\begin{split}
\|v_{k,r}\|^2_{H^1(B_{1/2})} \le \frac{1}{\big (\min(a,b) \big )^2}\|V_{k,r}\|^2_{H^1(B_{1/2})}\le c_1(N,M,\mu,r)\Big ( \|v\|^2_{H^1(B_{1/2})}+1 \Big ).
\end{split}
\ee
Since $u_k \wto u_0$ weakly in $H^1(B_{1/2})$, Lemma \ref{prop of Tab}-\ref{p2.1} gives $U_k\wto U_0$ weakly in $H^1(B_{1/2})$. From the properties of weak convergence in $H^1(B_{1/2})$ and Sobolev embeddings, we can verify that the following convergences hold up to a subsequence
\begin{enumerate}[label = \textbf{(C\arabic*)}, ref= C\arabic*]
\item \label{C0.5} $\nabla u_k \wto \nabla u_0$ and $\nabla U_k \wto \nabla U_0$ weakly in $H^1_0(B_{1/2})$ (c.f. Proposition \ref{prop of Tab}-\ref{p2.1}).
\item \label{C1} $u_k \to u_0$ and $U_k \to U_0$ strongly in $L^2(B_{1/2})$ (c.f. Proposition \ref{prop of Tab}-\ref{3b}).
\item \label{C2} $u_k\to u_0$ and $U_k \to U_0$ pointwise almost everywhere in $B_{1/2}$ (c.f. Proposition \ref{prop of Tab}-\ref{p7}).
\end{enumerate}
In order to proceed with the proof, we need some auxiliary estimates which are listed below
\begin{enumerate}[label=\textbf{(E\Roman*)}]
\item\label{E1}
For every $w\in H^1(B_{1/2})$ we have
\be\label{e1}
\GG_k(w;B_{1/2}) \le \frac{C(N,\lambda_+)}{k}\Big (\| w \|_{H^1(B_{1/2})}^2+1\Big ).
\ee
\item \label{E2} For $U_0 = \TT_{a,b}(u_0)$ we have 
\be\label{char u0} 
\int_{B_{1/2}} Q(0)\lambda(U_0)\,dx\le \liminf_{k\to \infty} \int_{B_{1/2}}Q(0)\lambda(U_k)\,dx.
\ee
\item \label{E3} For $V_{k,r} := \TT_{a,b}(v_{k,r})$,

\be\label{charv p=2}
\limsup_{r\to 0} \left ( \limsup_{k\to \infty} \int_{B_{1/2}} Q(0)\lambda(V_{r,k})\,dx \right )\le \int_{B_{1/2}}Q(0)\lambda(V)\,dx.
\ee
\item\label{E4} For $U_k := \TT_{a,b}(u_k)$ and $V:= \TT_{a,b}(v)$ we have for every $k\in \N$
\be\label{e4}
\begin{split}
\liminf_{k\to \infty} \int_{B_{1/2}}\big (|\nabla U_k|^2 - |\nabla V_{k,r}|^2 \big )\,dx  &\ge \int_{B_{1/2}}\big (|\nabla U_0|^2-|\nabla V|^2 \big )\,dx.
\end{split}
\ee

\end{enumerate}
Before delving into the proofs of \ref{E1}, \ref{E2}, \ref{E3} and \ref{E4}, we assume for the moment that all of them hold and prove minimality of $u_0$ for $J_0(\cdot;B_{1/2})$. We come back to their proofs after that.

First, we recall that for any $w\in H^1(B_{1/2})$ we have (c.f. \eqref{lambda preserved} and \eqref{Tab Lp nabla})
\be\label{equivalent}
\begin{split}
J_0(w;B_{1/2})  &= \int_{B_{1/2}} \big (a^2 |\nabla w^+|^2+ b^2|\nabla w^-|^2 + Q(0) \lambda(w) \big )\,dx\\
&= \int_{B_{1/2}}\left ( |\nabla \TT_{a,b} (w)| ^2 + Q(0) \lambda(\TT_{a,b}(w)) \right )\,dx.
\end{split}
\ee
From Lemma \ref{corollary to lemma 1.3}, we observe that to prove that $u_0$ is an absolute minimizer of $J_0(\cdot;B_{1/2})$, it is equivalent to show that $U_0 = \TT_{a,b}(u_0)$ is an absolute minimizer of $\FF_0(\cdot;B_{1/2})$ which is given by 
$$
\FF_0(W;B_{1/2}) := \int_{B_{1/2}}\big ( |\nabla W|^2+ Q(0)\lambda(W) \big )\,dx.
$$
So, our task is reduced to prove that $U_0$ is an absolute minimizer of $\FF_0(\cdot;B_{1/2})$. We make the following claim 
\be\label{M'}\tag{$\pmb{\MM'}$}
\mbox{\textbf{Claim:} $U_0$ is an absolute minimizer of $\FF_0(\cdot;B_{1/2})$.}
\ee
From \eqref{minimality of uk}, \eqref{min of uk 1} and \eqref{equivalent}, we have
$$
\int_{B_{1/2}}\Big ( |\nabla U_k|^2 + Q(0)\lambda(U_k) \Big )\,dx + \GG_k(u_k) \le \int_{B_{1/2}}\Big ( |\nabla V_{k,r}|^2 + Q(0)\lambda(V_{k,r}) \Big )\,dx + \GG_k(v_{k,r}).
$$
We rearrange the terms in the inequality above and take $\liminf_{k\to \infty}$ on the both sides to arrive to
\be\label{min of Uk}
\begin{split}
 \liminf_{k\to \infty} \Bigg [\int_{B_{1/2}}\Big (|\nabla U_k|^2 -|\nabla V_{k,r}|^2 \Big )\,dx + \int_{B_{1/2}} Q(0) \Big ( \lambda(U_k)- \lambda(V_{k,r})\Big )\,dx \Bigg ]  \le   \liminf_{k\to \infty}\big (\GG_{k}(v_{k,r}) -\GG_{k}(u_k) \big ). 
\end{split}
\ee
From \eqref{bound on uk}, \eqref{bound on vkr} and \ref{E1} we have 
\be\label{claim for vk p=2}
\lim_{k\to \infty} \GG_k(v_{k,r};B_{1/2}) = \lim_{k\to \infty} \GG_k(u_k;B_{1/2})=0.
\ee
Thus, from \eqref{min of Uk} we conclude
\be\label{ineq 1}
\liminf_{k\to \infty} \Bigg [\int_{B_{1/2}}\Big (|\nabla U_k|^2 -|\nabla V_{k,r}|^2 \Big )\,dx + \int_{B_{1/2}} Q(0) \Big ( \lambda(U_k)- \lambda(V_{k,r})\Big )\,dx \Bigg ]   \le 0.
\ee
Moreover, from \eqref{bounds on Uk and U0} and \eqref{bound on Vkr} the first integral in the estimates above is bounded independent of $k$, i.e. for every $k\in \N$ 
$$
\Big | \int_{B_{1/2}} \Big (|\nabla U_k|^2 -|\nabla V_{k,r}|^2 \Big )\,dx \Big | \le C(N,M, \mu, r)\Big ( \|v\|^2_{H^1(B_{1/2})}+1 \Big ).
$$
Furthermore,  
$$
 \Big | \int_{B_{1/2}} Q(0) ( \lambda(U_k)- \lambda(V_{k,r}))\,dx \Big | \le 2 C(N) q_2\lambda_{+}.
$$
{Therefore, we are entitled to use the \ref{I3}* in \eqref{ineq 1} to obtain }
$$
\liminf_{k\to \infty} \int_{B_{1/2}}\Big (|\nabla U_k|^2 -|\nabla V_{k,r}|^2 \Big )\,dx + \liminf_{k\to \infty} \int_{B_{1/2}} Q(0) \Big ( \lambda(U_k)- \lambda(V_{k,r}) \Big )\,dx   \le 0.
$$
Now rearranging the terms and {from \ref{I1}* and \ref{I2}**  together with \ref{E4} we have }
\be\label{ineq 2}
\begin{split}
 \int_{B_{1/2}}(|\nabla U_0|^2-|\nabla V|^2)\,dx &\le  \liminf_{k\to \infty} \int_{B_{1/2}} \Big (|\nabla U_k|^2 -|\nabla V_{k,r}|^2 \Big )\,dx\\
  & \le   \limsup_{k\to \infty}  \int_{B_{1/2}} Q(0) \Big (\lambda(V_{k,r})- \lambda(U_k) \Big )\,dx \\
  &\le \limsup_{k\to \infty} \int_{B_{1/2}}Q(0)\lambda(V_{k,r})\,dx + \limsup_{k\to \infty}\int_{B_{1/2}} -Q(0)\lambda(U_k)\,dx\\
&=   \limsup_{k\to \infty}  \int_{B_{1/2}}Q(0) (\lambda(V_{k,r})\,dx   -\liminf_{k\to \infty} \int_{B_{1/2}} Q(0)\lambda(U_k)\,dx.
\end{split}
\ee
Taking $\limsup_{r\to 0}$ on both sides of the inequality above and using \ref{E2} and \ref{E3}, we arrive at
\be\label{min of Uk 2}
\begin{split}
 \int_{B_{1/2}}(|\nabla U_0|^2-|\nabla V|^2)\,dx & \le  \limsup_{r\to 0}  \Bigg (  \limsup_{k\to \infty}  \int_{B_{1/2}}Q(0) \lambda(V_{k,r})\,dx   -\liminf_{k\to \infty} \int_{B_{1/2}} Q(0)\lambda(U_k)\,dx \Bigg )\\
 &\le \limsup_{r\to 0} \Bigg (\limsup_{k\to \infty}  \int_{B_{1/2}}Q(0) \lambda(V_{k,r})\,dx\Bigg ) - \liminf_{k\to \infty} \int_{B_{1/2}}Q(0) \lambda(U_k)\,dx\\
&\le   \int_{B_{1/2}} Q(0)\Big (\lambda(V) -\lambda(U_0) \Big )\,dx.
\end{split}
\ee
Hence, we find 
$$
\int_{B_{1/2}}(|\nabla U_0|^2-|\nabla V|^2)\,dx\le \int_{B_{1/2}} Q(0)(\lambda(V) -\lambda(U_0))\,dx.
$$
This proves $U_0$ is an absolute minimizer of $\FF_0(\cdot; B_{1/2})$,  i.e. Claim \ref{M'} is proven.
In order to finish the proof of Theorem \ref{step 1}, Now, let us prove \ref{E1}-\ref{E4}.  To prove \ref{E1}, we observe that since $\Linfty{A_{\pm}-A_{\pm}(0)}{B_{1/2}}\le \frac{1}{k}$, then
\be\label{part 1}
\begin{split}
\int_{B_{1/2}}\Big ((A_{\pm}(0) -A_{\pm,k}) |\nabla v|^2 \Big ) \,dx &\le \Linfty{A_{\pm}(0) -A_{\pm,k}}{B_{1/2}} \int_{B_{1/2}}|\nabla v|^2 \,dx\\
& \le \frac{\L{2}{\nabla v}{B_{1/2}}^2}{k} \le \frac{\|v\|_{H^1(B_{1/2})}^2}{k}.
\end{split}
\ee 
Moreover, since $\Linfty{Q_k -Q(0)}{B_{1/2}} \le\frac{1}{k} $, we have 
\be\label{part 2}
\int_{B_{1/2}} (Q_k -Q(0)) \lambda(v)\,dx \le C(N)\lambda_+ \Linfty{Q_k -Q(0)}{B_{1/2}}  \le \frac{C(N)\lambda_+}{k}.
\ee
Let $F_k:= |f_{+,k}| + |f_{-,k}|$. In the case $N>2$, we have ${2^*}' =\frac{2N}{N+2}<N$. So, by Sobolev embedding, Young's inequality and \eqref{def of coeff} we find
\be\label{fpmv p=2}
\begin{split}
\int_{B_{1/2}} |f_k(x,v_k)v|\,dx &\le \int_{B_{1/2}} \Big (|f_{+,k}| + |f_{-,k}| \Big )|v| \,dx = \int_{B_{1/2}} |F_k| |v|\,dx \\
&\le \|F_k\|_{L^{{2^*}'}(B_{1/2})} \|v\|_{L^{2^*}(B_{1/2})}\\
&\le  \frac{C(N)}{k}\|v\|_{H^1(B_{1/2})}\le \frac{C(N)}{2k}\Big (\|v\|_{H^1(B_{1/2}  )}^2+1\Big ).\\
\end{split}
\ee 
If $N=2$, we know that $H^1(B_{1/2})$ embeds in $L^{N'}(B_{1/2})$ and therefore
\be\label{fpmv p=2 N=2}
\begin{split}
\int_{B_{1/2}} |f_k(x,v_k)v|\,dx &\le \int_{B_{1/2}} (|f_{+,k}| + |f_{-,k}| )|v| \,dx = \int_{B_{1/2}} |F_k| |v|\,dx \\
&\le \|F_k\|_{L^{{N}}(B_{1/2})} \|v\|_{L^{N'}(B_{1/2})}\\
&\le  \frac{C(N)}{k}\|v\|_{H^1(B_{1/2})}\le \frac{C(N)}{2k}(\|v\|_{H^1(B_{1/2})}^2+1).\\
\end{split}
\ee

{From \eqref{part 1}, \eqref{part 2} and \eqref{fpmv p=2} or \eqref{fpmv p=2 N=2}} the claim \ref{E1} is proven. 
To prove \ref{E2} we observe that 
\be\label{claim 2v p=2}
\lambda(U_0)=\lambda_+ \chi_{\{U_0>0\}}+ \lambda_- \chi_{\{U_0\le 0\}} \le \liminf_{k\to \infty}\Big (\lambda_+ \chi_{\{U_k>0\}}+ \lambda_- \chi_{\{U_k\le 0\}} \Big ) = \liminf_{k\to \infty}\lambda(U_k)\mbox{ a.e. in $B_{1/2}$}
\ee
Indeed, since $\lambda_- <\lambda_+$, we have $\lambda_- \le \liminf_{k\to \infty}( \lambda_+ \chi_{\{U_k>0\}}+ \lambda_- \chi_{\{U_k\le 0\}})$. Hence \eqref{claim 2v p=2} holds in the set $\{U_0=0\}$.  From \eqref{C2} we have $U_k\to U_0$ pointwise almost everywhere in $B_{1/2}$. Suppose $x_0 \in \{U_0>0\}$ is such that $U_k(x_0) \to U_0(x_0)$. Then, for $k$ sufficiently large we have $x_0 \in \{U_k>0\}$. Hence \eqref{claim 2v p=2} holds almost everywhere in $\{U_0>0\}$. Similarly we can show the same for the set $\{U_0<0\}$. This proves the claim (\ref{claim 2v p=2}).  Since $Q(0)\lambda(U_k) \le q_2\lambda_+$ for all $k$. By \eqref{claim 2v p=2} and since $Q(0)\ge 0$ we have by Fatou's lemma  
\[
\begin{split}
\int_{B_{1/2}}Q(0)\lambda(U_0)\,dx&\le \int_{B_{1/2}} Q(0)\liminf_{k\to \infty}\lambda(U_k)\,dx \le \liminf_{k\to \infty} \int_{B_{1/2}} Q(0)\lambda(U_k)\,dx.\\
\end{split}
\]
This proves \ref{E2}. Now we prove \ref{E3}. The integral in the LHS of \ref{E3} can be written as 
$$
\int_{B_{1/2}} Q(0)\lambda(V_{r,k})\,dx =Q(0) \Big ( \lambda_+|\{V_{k,r}>0\} \cap B_{1/2}| + \lambda_-|\{V_{k,r}\le 0\} \cap B_{1/2}| \Big ).
$$
In order to prove \ref{E3}, we first observe 
\be\label{charv2}
\begin{split}
&\lim_{r\to 0} \left ( \limsup_{k\to \infty} |\{V_{k,r}>0\} \cap B_{1/2}| \right ) = |\{V>0\} \cap B_{1/2}|\;\;\mbox{and}\\
&\lim_{r\to 0} \left ( \limsup_{k\to \infty} |\{V_{k,r}\le 0\} \cap B_{1/2}| \right ) = |\{V\le 0\} \cap B_{1/2}|.
\end{split}
\ee
Indeed, let us set $\RR_{r}:= B_{1/2}\setminus B_{\frac{1}{2}-r}$. Now, since $\eta_r=1$ in $B_{1/2}\setminus \RR_r = B_{\frac{1}{2}-r}$ (c.f. \eqref{eta}), recalling the definition of $V_{k,r}$ in \eqref{def of  Vkr} we observe that
$$
V_{k,r} = V \;\mbox{ in $B_{1/2}\setminus \RR_r$.} 
$$
With this in mind we have
\be\label{calc for vkr}
\begin{split}
|\{V_{r,k}>0 \}\cap B_{1/2}| &= |( \{V_{r,k}>0 \} \cap B_{1/2})\setminus \RR_r |+ |\{V_{r,k}>0 \}\cap \RR_{r}|\\
&= |(\{V>0 \}\cap B_{1/2})\setminus \RR_r |+ |\{V_{r,k}>0 \}\cap \RR_{r}|\\
&= |\{V>0 \}\cap B_{1/2}| - | \RR_{r}\cap \{V>0\}| + |\{V_{r,k}>0 \}\cap \RR_{r}|.
\end{split}
\ee
From above computations we have
\[
\begin{split}
|\{V>0\}\cap B_{1/2}| - |\RR_{r}\cap \{V>0\}| \le |\{V_{r,k}>0\}\cap B_{1/2}| \le |\{V>0\}\cap B_{1/2}| + |\{V_{r,k}>0 \}\cap \RR_{r}|.
\end{split}
\]
and therefore
$$
|\{V>0\}\cap B_{1/2}| - |\RR_{r}| \le |\{V_{r,k}>0\}\cap B_{1/2}| \le |\{V>0\}\cap B_{1/2}| + | \RR_{r}|.
$$
Now, taking $\limsup_{k\to \infty}$ and $\limsup_{r\to 0}$ in the inequality above we find 
\[
\begin{split}
|\{V>0\}\cap B_{1/2}| &= \lim_{r\to 0}  \Big (|\{V>0\}\cap B_{1/2}| - |\RR_{r}| \Big )\\
 &\le \limsup_{r\to 0} \Big ( \limsup_{k\to \infty} |\{V_{k,r}>0\}\cap B_{1/2}| \Big  )\\
& \le \limsup_{r\to 0} \Big (|\{V>0\}\cap B_{1/2}| + | \RR_{r}| \Big ) =|\{V>0\}\cap B_{1/2}|.
\end{split}
\]
We can proceed similarly by replacing $\{V_{k,r}>0\}$ by $\{V_{k,r}\le 0\}$ in \eqref{calc for vkr} to arrive to 
$$
\limsup_{r\to 0} \Big ( \limsup_{k\to \infty} |\{V_{k,r}\le 0\}\cap B_{1/2}| \Big  ) =|\{V\le0\}\cap B_{1/2}|.
$$
This proves \eqref{charv2}. {Since the sequences $|\{V_{k,r}>0\}|$ and $|\{V_{k,r}\le 0\}|$ are bounded in $k$ and $r$, we can use \ref{I4}** twice (first in $k$ then in $r$)} and \eqref{charv2} to obtain
\[
\begin{split}
\limsup_{r\to 0} &\left ( \limsup_{k\to \infty} \int_{B_{1/2}} Q(0)\lambda(V_{r,k})\,dx \right )\\
 &= Q(0) \limsup_{r\to 0} \Bigg ( \limsup_{k\to \infty} \Big ( \lambda_+ |\{V_{k,r}>0\}\cap B_{1/2}| + \lambda_- |\{V_{k,r}\le 0\}\cap B_{1/2}|\Big ) \Bigg )\\
&\le Q(0) \limsup_{r\to 0} \Bigg (  \limsup_{k\to \infty} \lambda_+ |\{V_{k,r}>0\}\cap B_{1/2}|  \Bigg )+ Q(0)\limsup_{r\to 0} \left (  \limsup_{k\to \infty} \lambda_- |\{V_{k,r}\le 0\}\cap B_{1/2}|  \right )\\
&=Q(0) \big (  \lambda_+ |\{V>0\}\cap B_{1/2}| + \lambda_- |\{V\le 0\}\cap B_{1/2}| \big ) \\
&= \int_{B_{1/2}}Q(0)\lambda(v)\,dx.
\end{split}
\]
This shows \ref{E3}. Finally, we prove \ref{E4}. We start by observing that 
\[
\begin{split}
\int_{B_{1/2}}\big (|\nabla U_k|^2 - |\nabla V_{k,r}|^2 \big )\,dx&= \int_{B_{1/2}} \big ( \nabla (U_k -V_{k,r}) \cdot \nabla (U_k+V_{k,r}) \big )\,dx.
\end{split}
\]
Now we plug in the definition of $V_{k,r}$ (c.f. \eqref{def of Vkr}) into the previous equation and arrive to 
\be\label{a+b}
\begin{split}
\nabla (U_k - V_{k,r}) &= \nabla (U_k - V -(1-\eta_r)(U_k-U_0))\\
&= \nabla (U_0-V +\eta_r(U_k-U_0))
\end{split}
\ee
and 
\be\label{a-b}
\begin{split}
\nabla (U_k+V_{k,r}) & = \nabla (U_k + V +(1-\eta_r)(U_k-U_0))\\
&= \nabla (U_0+V + (2-\eta_r)(U_k-U_0)).
\end{split}
\ee
From \eqref{a+b} and \eqref{a-b} we write 
\be\label{all the terms}
\begin{split}
\int_{B_{1/2}} \nabla (U_k -V_{k,r}) \cdot \nabla (U_k+V_{k,r})\,dx = \int_{B_{1/2}}  \nabla &(U_0-V)\cdot \nabla (U_0+V)\,dx+\mathbf{T^1_{k}}+\mathbf{T_{k}^2}+\mathbf{T^3_{k}}.
\end{split}
\ee
where
\be\label{the terms}
\begin{split}
\mathbf{T^1_{k}} &:= \int_{B_{1/2}}\nabla (U_0 -V) \cdot \nabla ((2-\eta_r)(U_k-U_0))\,dx \\
\mathbf{T^2_{k}}&:= \int_{B_{1/2}}\nabla(\eta_r(U_k-U_0)) \cdot \nabla (U_0+V)\,dx \\
\mathbf{T^3_{k}}&:= \int_{B_{1/2}}\nabla ( \eta_r(U_k-U_0))\cdot \nabla ((2-\eta_r)(U_k-U_0))\,dx.
\end{split}
\ee
Let us study each term separately. First, we see that
\[
\begin{split}
\mathbf{T^1_{k}} =\int_{B_{1/2}} \nabla (U_0 -V) \cdot \nabla ((2-\eta_r)(U_k-U_0))\,dx&= \int_{B_{1/2}} (2-\eta_r) \nabla (U_0-V)\cdot \nabla (U_k-U_0)\,dx\\
&-\int_{B_{1/2}} \big ( (U_k -U_0) \nabla (U_0-V) \cdot \nabla \eta_r \big )\,dx.
\end{split}
\]
In the identity above, the first term on the right hand side tends to zero as $k\to \infty$ since from \eqref{C0.5}, $\nabla (U_k - U_0)$ goes to zero weakly in $L^2(B_{1/2})$. Also, the second term on the right hand side tends to zero as $k\to \infty$ because from \eqref{C1}, $U_k$ converges to $U_0$ strongly in $L^2(B_1)$. Therefore 
\be\label{14}
\lim_{k\to \infty} \mathbf{T^1_{k}}= \lim_{k\to \infty}\int_{B_{1/2}} \nabla (U_0 -V) \cdot \nabla ((2-\eta_r)(U_k-U_0))\,dx =0.
\ee
Now, let us look at the following term
\[
\begin{split}
\mathbf{T^2_{k}}= \int_{B_{1/2}} \nabla( \eta_r(U_k-U_0)) \cdot \nabla (U_0+V) \,dx = \int_{B_{1/2}} &\eta_r \nabla (U_0+V)\cdot \nabla (U_k-U_0)\,dx\\
& + \int_{B_{1/2}} (U_k-U_0)\nabla \eta_r\cdot \nabla (U_0+V)\,dx.
\end{split}
\]
Thus, by similar arguments as those in \eqref{14}, we obtain that 
\be\label{23}
\lim_{k\to \infty} \mathbf{T^2_{k}} = \lim_{k\to \infty} \int_{B_{1/2}} \nabla( \eta_r(U_k-U_0)) \cdot \nabla (U_0+V) \,dx =0.
\ee
Now, we estimate $\mathbf{T^3_{k}}$. We write $\mathbf{T^3_{k}}$ as follows 
\be\label{*1}
\begin{split}
\mathbf{T^3_{k}} = \mathbf{S^3_{k}}+\mathbf{T^3_{k}}^*.
\end{split}
\ee
where 
\[
\begin{split}
\mathbf{S^3_{k}}:=  \int_{B_{1/2}}\eta_r(2-\eta_r) |\nabla (U_k-U_0)|^2\,dx 
\end{split}
\]
and 
$$
\mathbf{T^3_{k}}^* := \int_{B_{1/2}} \Big ( (2-\eta_r)(U_k-U_0) \nabla (U_k-U_0) \cdot \nabla \eta_r- \eta_r (U_k-U_0)\nabla (U_k-U_0)\cdot \nabla \eta_r- |U_k-U_0|^2 |\nabla \eta_r|^2 \Big ) \,dx.
$$
From \eqref{C1} and \eqref{C2} we have 
\be\label{*2}
\lim_{k\to \infty}\mathbf{T^3_{k}}^* =0.
\ee
Since $0\le \eta_r\le 1$ we have for every $k\in \N$
\be\label{*3}
\mathbf{S^3_{k}} = \int_{B_{1/2}} \Big ( \eta_r(2-\eta_r) |\nabla (U_k-U_0)|^2\Big )\,dx \ge 0.
\ee
Therefore from \eqref{*3} $\mathbf{T^3_k}\ge \mathbf {{T^3_k}^*}$ for all $k$. From \eqref{*1} and \eqref{*2} we obtain
\be\label{24}
\liminf_{k\to \infty} \mathbf{T^3_{k}}  \ge \liminf_{k\to \infty} \mathbf{T^3_{k}}^*=0.
\ee
Now, by using \eqref{14}, \eqref{23}, \eqref{24} and obtain 
\[
\begin{split}
\liminf_{k\to \infty} \int_{B_{1/2}}\big (|\nabla U_k|^2 - |\nabla V_{k,r}|^2 \big )\,dx 
&= \liminf_{k\to \infty} \Big ( \int_{B_{1/2}} \big (|\nabla U_0|^2 -|\nabla V|^2 \big )\,dx +  \mathbf{T^1_{k}}+\mathbf{T_{k}^2}+\mathbf{T^3_{k}} \Big )\\
&\ge \int_{B_{1/2}} (U_0-V)\cdot \nabla (U_0+V) \,dx= \int_{B_{1/2}}(|\nabla U_0|^2-|\nabla V|^2)\,dx.
\end{split}
\]
This shows \ref{E4} and this finally proves $u_0$ is an absolute minimizer of $J_0(\cdot;B_{1/2})$. This is a contradiction to the assumption (\ref{absurdum}). Thus, we finish the proof of Theorem \ref{step 1}.
}
\end{proof}

Now, we prove that the value of $\delta>0$ in the Lemma \ref{step 1} is dependent only on $L^{\infty}$ bounds of absolute minimizer $u$ and other  universal constants related to the problem. For this we use the Widman's hole filling technique by following ideas from \cite[Chapter 7]{eg05} (see also \cite{W71}).

\begin{proposition}\label{step 1.5}
Let $u\in H^1(B_1)\cap L^{\infty}(B_1)$ be an absolute minimizer of $J(\cdot; B_1)$ such that 

\noindent $\Linfty{u}{B_1}\le 1$. Then for every $\eps>0$ there exists $0<\delta(\eps, \mu, q_2,\lambda_+, N)<1$ such that if  
$$
 \max \Big (\|A_{\pm}-A_{\pm}(0)\|_{L^{\infty}(B_1)},\, \|f_{\pm}\|_{L^N(B_1)},\, \|Q-Q(0)\|_{L^{\infty}(B_1)} \Big ) \le \delta 
$$
then 
$$
\|u-h\|_{L^{\infty}(B_{1/4})}\le\eps.
$$
Here $h\in H^{1}(B_{1/2})\cap L^{\infty}(B_{1/2})$ is such that $\Linfty{h}{B_{1/2}}\le 1$ and $h$ is also an absolute minimizer of $$J_0(v;B_{1/2}):= \int_{B_{1/2}}\Big (A(0,v)|\nabla v|^2 + Q(0)\lambda(v) \Big )\,dx.$$ 
\end{proposition}
\begin{proof}
We observe that $|f(x,s)| \le F:=(|f_{+}| + |f_{-}|)$ and $|Q(x)\lambda(s)|\le q_2\lambda_+$. Hence, we can easily verify that the integrand of the functional $J(\cdot;B_1)$ satisfies the following for almost every $x\in B_1$ and for all $\xi\in \R^N$, $s\in \R$.
\be\label{lbub}
{\mu}|\xi|^2 -F |s| - q_2\lambda_+ \le A(x,s)|\xi|^2 +f(x,s)s+Q(x)\lambda(s) \le \frac{1}{\mu} |\xi|^2 + F |s| + q_2\lambda_+
\ee
where $F:= |f_{+}| + |f_{-}|$.  Let $r,s>0$ be such that $1/2\le r<s\le \frac{3}{4}$. Since, $u$ is an absolute minimizer of $J(\cdot;B_1)$, so it is of $J(\cdot;B_s)$ (c.f. Lemma \ref{local minimizer}). This means that for every $v\in H^1(B_s)$ such that $v-u\in H_0^1(B_s)$, we have $J(u;B_s)\le J(v;B_s)$. In other words
\be\label{minimality in Bs}
\int_{B_s}\Big ( A(x,u)|\nabla u|^2 +f(x,u)u+Q(x)\lambda(u) \Big )\,dx \,\le \, \int_{B_s} \Big (A(x,v)|\nabla v|^2 +f(x,v)v+Q(x)\lambda(v) \Big )\,dx.
\ee
Therefore, from \eqref{lbub} we can write 
\be\label{min of uk}
\int_{B_s} \left ( {\mu}|\nabla u|^2 -F |u| -q_2\lambda_+  \right )\,dx\le \int_{B_s} \left ( \frac{1}{\mu}|\nabla v|^2 +F |v| +q_2\lambda_+  \right )\,dx
\ee
for every $v\in H^1(B_s)$ such that $v-u \in H_0^{1}(B_s)$. Let $\eta \in C_0^{\infty}(B_1)$ be such that 
$$
0\le \eta\le 1,\qquad \eta(x)=
\begin{cases}
1 \;\;x\in B_r\\
0\;\;x\in B_1\setminus B_s.
\end{cases}\qquad \mbox{\big (recall $\frac{1}{2}\le r<s\le \frac{3}{4}$ \big )}
$$
We can assume that $|\nabla \eta|\le \frac{C(N)}{|s-r|}$ in $B_1$.  Now, we consider the test function $v=u(1-\eta)$. Since, $v-u= -\eta u\in H_0^1(B_s)$, the function $v\in H^1(B_s)$ is an admissible function for the inequality \eqref{min of uk}.

Since $\LN{F}{B_s}\le 2\delta <2$ and $\L{N'}{\big ( |u|+|u(1-\eta)|\big )}{B_s}\le 2C(N)s^{\frac{N}{N'}}=2C(N)s^{N-1}$. By plugging $v=u(1-\eta)$ in to \eqref{min of uk} and rearranging the terms accordingly we obtain,
\[
\begin{split}
\int_{B_s} {\mu}|\nabla u|^2\,dx &\le \int_{B_s}\left ( \frac{1}{\mu}|\nabla (u(1-\eta))|^2 +F |u(1-\eta)| +q_2\lambda_+  \right )\,dx+ \int_{B_s}\left (F |u| +q_2\lambda_+  \right )\,dx\\
&\le \int_{B_s}\frac{1}{\mu} \Big |  (1-\eta)\nabla u -u \nabla \eta  \Big |^2\,dx + \LN{F}{B_s} \L{N'}{\big ( |u|+|u(1-\eta)|\big )}{B_s}+C(N)q_2\lambda_+ s^N\\
&\le \frac{2}{\mu} \int_{B_s} \Big ( (1-\eta)^2|\nabla u|^2+|\nabla \eta|^2 |u|^2 \Big )\,dx + 4C(N) s^{N-1} + C(N,q_2,\lambda_+)s^N\\
&\le C(\mu)\int_{B_s}\Big ( (1-\eta)^2 |\nabla u|^2+ |\nabla \eta|^2|u|^2\Big )\,dx +C(N,q_2,\lambda_+)s^{N-1}. \\
\end{split}
\]
Since $\Linfty{u}{B_1}\le 1$, $s\in (0,1)$, $|\nabla \eta|\le \frac{C(N)}{|s-r|}$ in $B_s$ and $\eta=1$ in $B_r$ we have
\[
\begin{split}
\int_{B_s}|\nabla u|^2\,dx&\le \bar C(\mu)\int_{B_s}\Big ((1-\eta)^2 |\nabla u|^2 + |\nabla \eta|^2u^2 \Big )\,dx + C( N,\mu,q_2,\lambda_+) \\
&\le \bar C(\mu) \Bigg [\int_{B_s\setminus B_r}|\nabla u|^2\,dx + \frac{C(N)}{|s-r|^2}\int_{B_s} |u|^2\,dx\Bigg ] +C( N,\mu,q_2,\lambda_+) \\
&\le \bar C(\mu) \Bigg [\int_{B_s\setminus B_r}|\nabla u|^2\,dx + \frac{C(N)}{|s-r|^2}s^N\Bigg ] +C( N,\mu,q_2,\lambda_+).
\end{split}
\]
Since $r<s$ and $s\in (0,1)$ we can write 
$$
\int_{B_r}|\nabla u|^2\,dx\le \int_{B_s}|\nabla u|^2\,dx\le \bar C(\mu) \int_{B_s\setminus B_r}|\nabla u|^2\,dx + \frac{C(N,\mu)}{|s-r|^2} +C(N,\mu,q_2,\lambda_+) .
$$ 
We now add $\bar C(\mu) \int_{B_r}|\nabla u|^2$ on LHS and RHS sides of the inequality above to obtain 
$$
\int_{B_r}|\nabla u|^2\,dx \le \frac{\bar C(\mu)}{\bar C(\mu)+1} \int_{B_s}|\nabla u|^2\,dx + \frac{C'(N,\mu)}{|s-r|^2}+C'(N,\mu,q_2,\lambda_+).
$$
Now, from Lemma \ref{useful lemma} with $Z(t):= \int_{B_t}|\nabla u|^2\,dx$ we conclude that 
\be\label{energy estimates} 
\int_{B_{1/2}}|\nabla u|^2 \le C_2(N, \mu,q_2,\lambda_+).
\ee
Once we have a universal estimate for the energy in $B_{1/2}$, we can apply Theorem \ref{step 1}. Since $u$ is an absolute minimizer of $J(\cdot;B_1)$, so it is of $J(\cdot;B_{1/2})$. Moreover, $\Linfty{u}{B_{1/2}} \le \Linfty{u}{B_1}\le 1$. Now, Theorem \ref{step 1} with $M = C_2(\mu,q_2,\lambda_+,N)^{1/2}>0$ provides a $\delta:=\delta(\eps,\mu,q_2,\lambda_+,N)>0$ for which Proposition \ref{step 2} holds. This finishes the proof.
\end{proof}

\section{Optimal regularity of minimizers}\label{main result}
We recall that in this section $J_{A,f,Q}(\cdot;B_1)$ satisfy the structural conditions \ref{G1}-\ref{G5}.
\begin{remark}\label{rescaling}
Let $u\in H^1(B_\Theta(x_0))$  be an absolute minimizer of $J(\cdot; B_\Theta(x_0))$. We define $w$ as follows 
$$
w(x) := \Phi u(x_0+\Theta x ).
$$
We can verify that $w$ is an absolute minimizer the following functional
$$
\bar J (v):= \int_{B_{r / \Theta}}\left (  \bar A(x,v) |\nabla  v|^2 -\bar f(x,v )v+\bar Q(x)\lambda (v)\right )\,dx
$$
where $\bar A$, $\bar f$ and $\bar Q$ are defined in $B_{r/\Theta}$ as follows
\[
\begin{split}
\bar A(x,s)&= A(x_0+\Theta x, s)\\
\bar f(x,s)&= \Phi\Theta^2 f(x_0+\Theta x)\\
\bar Q(x)&= \Phi^2 \Theta^2 Q(x_0+\Theta x).
\end{split}
\]
\end{remark}
We now prove the key lemma analogous to Proposition \ref{step 0.1} in the context of minimizers of $J_{A,f,Q}(\cdot;B_1)$.

\begin{proposition}[Key lemma for minimizers]\label{step 2}
Let $u\in H^1(B_1)\cap L^\infty(B_1)$ be an absolute minimizer of $J(\cdot; B_1)$ satisfying the structural condition \ref{G1}-\ref{G5} with $u(0)=0$ and 

\noindent $\Linfty{u}{B_1}\le 1$.  Then for any $0<\alpha<1$, there exists $\delta (N,\mu, q_2,\lambda_+,\alpha) >0$ and $0< R_0(N,\mu, q_2,\lambda_+,\alpha)<1/4$ such that if 
$$
 \max (\|A_{\pm}-A_{\pm}(0)\|_{L^{\infty}(B_1)},\, \|f_{\pm}\|_{L^N(B_1)},\, \|Q-Q(0)\|_{L^{\infty}(B_1)}) \le \delta 
$$
then we have 
\be\label{claim step 2}
\sup_{B_{R_0}}|u|\le R_0^{\alpha}.
\ee
\end{proposition}
\begin{proof}
Let $\eps>0$ which will be suitably chosen later. We know that for $\delta(\eps)>0$ and $h$ as in Lemma \ref{step 1} we have  
\be\label{s2e1}
\|u-h\|_{L^{\infty}(B_{1/2})}<\eps.
\ee 
Fix $\beta=\frac{1+\alpha}{2}$ from Proposition \ref{lemma 1.3} we have
\be\label{s2e1.}
\sup_{B_r}|h-h(0)|\le C(N,\mu,q_2,\lambda_+,\alpha)r^{\beta}, \;\;r<\frac{1}{4}.
\ee

By (\ref{s2e1}) and (\ref{s2e1.}) we get for all $r<1/4$
\be\label{s2e3}
\begin{split}
\sup_{B_r}|u(x)-u(0)|&\le \sup_{B_r}\Big ( |u(x)-h(x)|+|h(x)-h(0)|+|h(0)-u(0)| \Big ) \\
&\le 2\eps+C(N, \mu, q_2,\lambda_+,\alpha)r^{\beta}.
\end{split}
\ee
We can select $R_0(N, \mu, q_2,\lambda_+,\alpha)<1/4$ such that 
$$
C(N, \mu, q_2\lambda_+,\alpha)R_0^{\beta}\le\frac{R_0^{\alpha}}{3}.
$$
This means that 
$$
R_0\le\Big (  \frac{1}{3C} \Big )^{2/(1-\alpha)}.
$$
Now, we can make our choice of $\eps$. Choose $\eps(N, \mu, q_2,\lambda_+,\alpha)>0$ in such a way that 
$$
\eps<\frac{R_0^{\alpha}}{3}.
$$
We see that the choice of $\delta$ depends on $\eps$ and since $\eps$ depends on $N, \mu, q_2,\lambda_+$ and $\alpha$ so does $\delta$. Since $u(0)=0$ and $C(N,\mu, q_2,\lambda_+,\alpha)R_0^{\beta}$ and $\eps$ are bounded by $R_0^{\alpha}/3$. We obtain from \eqref{s2e3}
$$
\sup_{B_{R_0}} |u| \le R_0^{\alpha}.
$$
\end{proof} 
We now have the ingredients to show the $C^{0,1^-}$ regularity estimates for an absolute minimizer $u$ around the points in the zero level set.
\begin{proposition}\label{step 3}
Let $u\in H^1(B_1)\cap L^{\infty}(B_1)$ be an absolute minimizer of $J(\cdot; B_1)$ with $u(0)=0$ and $\Linfty{u}{B_1}\le 1$.  Then for every $0<\alpha<1$, there exists a $\delta(N, \mu, q_2,\lambda_+,\alpha)>0$ and $C(N, \mu, q_2,\lambda_+,\alpha)>0$ such that if 
$$
 \max \Big (\|A_{\pm}-A_{\pm}(0)\|_{L^{\infty}(B_1)},\, \|f_{\pm}\|_{L^N(B_1)},\, \|Q-Q(0)\|_{L^{\infty}(B_1)} \Big ) \le \delta 
$$
then for $R_0(N,\mu, q_2,\lambda_+,\alpha)$ as in Proposition \ref{step 2} we have
\be\label{claim step 3}
\sup_{B_{r}}|u(x)|\le C(N, \mu, q_2,\lambda_+,\alpha) \cdot r^{\alpha} \;\; \forall \, r\le R_0.
\ee
Precisely speaking, we have $C(N, \mu, q_2,\lambda_+,\alpha)= R_0^{-\alpha}$.
\end{proposition}

\begin{proof}
We argue by scaling and claim that for all $k\in \N$
\be\label{s3e1}
 \sup_{B_{R_0^{k}}}|u(x)|\le R_0^{k\alpha }.
\ee
It follows readily from Proposition \ref{step 2} that (\ref{s3e1}) holds for $k=1$. Now, let us assume it holds up to $k_0\in \N$. We prove that (\ref{s3e1}) also holds for $k=k_0+1$. In order to do that, we set the following scaled function 
$$
\tilde u(y)=\frac{1}{R_0^{k_0\alpha}}u(R_0^{k_0}y).
$$
By Remark \ref{rescaling} we see that $\tilde u$ minimizes the functional $\tilde J$ given by
$$
\tilde J(v) :=  \int_{B_1} \tilde A(y,v)|\nabla v|^2 - \tilde f(y,v)v + \tilde Q(y)\lambda(v)\,dx
$$
where $\tilde A$, $\tilde f$ and $\tilde Q$ are defined as follows for $y\in B_1$
\[
\begin{split}
&\tilde A(y,s):=A(R_0^{k_0}y,s)\\
&\tilde f(y,s):=R_0^{(2k_0(1-\alpha)+k_0\alpha)}f(R_0^{k_0}y,s)\\
&\tilde Q (y):=R_0^{2k_0(1-\alpha)} Q (r_0^{k_0}y).
\end{split}
\]
We observe that $\tilde J$ and $\tilde u$ satisfy the assumptions of Proposition \ref{step 2}. Indeed, by induction hypothesis (from \eqref{s3e1} for $k=k_0$) we have 
$$
\sup _{B_{1}}|\tilde u|= R_0^{-k_0 \alpha} \sup_{B_{R_0^{k_0}}}|u|\le 1.
$$
Also, for $\delta>0$ as in Proposition \ref{step 2}, we can see that
$$
\sup_{B_1}|\tilde A_{\pm} - \tilde A_{\pm}(0)|=  \sup_{B_{R_0^{k_0}}} |A_{\pm}  -A_{\pm}(0)| \le \delta
$$
and 
$$
\sup_{B_1}|\tilde Q - \tilde Q(0)| =  R_0^{2k_0(1-\alpha)} \sup_{B_{R_0^{k_0}}} |Q  -Q(0)| \le \delta,\;\;\;\;|\tilde Q|\le q_2.
$$
Furthermore, 
$$
\|\tilde f\|_{L^{N}(B_1)}=R_0^{k_0(1-\alpha)}\|f\|_{L^N(B_{R_0^k})}\le  \delta.
$$
Moreover $\tilde u(0)=0$. By applying Proposition \ref{step 2} to the pair ($\tilde J$, $\tilde u$), we see from \eqref{claim step 2} that 
$$
\sup_{B_{R_0}}|\tilde u|\le R_0^{\alpha}.
$$
Rescaling back to $u$, we obtain 
$$
\sup_{B_{R_0^{k_0+1}}}|u|\le R_0^{(k_0+1)\alpha}.
$$
This finishes the proof by induction and proves the claim (\ref{s3e1}) for all $k\in \N$. To prove Proposition \ref {step 3} let us take $r\in (0,R_0)$ and $k$ such that $R_0^{k+1}\le r<R_0^k$. From (\ref{s3e1}), we see that 
$$
\sup_{B_r}|u|\le \sup_{B_{R_0^{k}}}|u|\le R_0^{k\alpha}=R_0^{(k+1)\alpha}\frac{1}{R_0^{\alpha}}\le \frac{1}{R_0^{\alpha}}r^{\alpha}.
$$
Now, we just take $C(N, \mu, q_2,\lambda_+,\alpha) = \frac{1}{R_0^{\alpha}}$. This finishes the proof of Proposition \ref{step 3}.
\end{proof}
Now we claim that only the smallness in oscillations of coefficients $A_{\pm}$ is sufficient to show the regularity estimates for an absolute minimizer $u$ of $J(\cdot;B_1)$. We prove this result in the following rescaled version of previous lemma.

\begin{proposition}\label{step 3.5}
Suppose $u\in H^1(B_{\rho}(x_0))\cap L^{\infty}(B_{\rho}(x_0))$ is an absolute minimizer of $J(\cdot; B_\rho(x_0))$ $(\rho<1)$ and $u(x_0)=0$.  Then for all $0<\alpha<1$, there exists $C(N, \mu, q_2,\lambda_+,\alpha)>0$ such that if 
\be\label{SMALLNESS}
\|A_{\pm}-A_{\pm}(x_0)\|_{L^{\infty}(B_\rho(x_0))} \le \delta 
\ee
for $\delta(N, \mu, q_2,\lambda_+,\alpha)>0$, then for $R_0(N, \mu, q_2,\lambda_+,\alpha)$ as in Proposition \ref{step 3} we have
\be\label{claim step 3}
\sup_{B_{r}(x_0)}|u(x)|\le \frac{C(N, \mu, q_2,\lambda_+,\alpha)}{\rho^{\alpha}} \left (  \rho+\Linfty{u}{B_{\rho}(x_0)} +\rho \LN{F}{B_\rho(x_0)} \right )  r^{\alpha} \;\; \forall\, r\le\rho R_0.
\ee
\end{proposition}
\begin{remark}
In fact, $C(N, \mu, q_2,\lambda_+,\alpha)$ in Proposition \ref{step 3.5} can be taken as 
$$
C(N, \mu, q_2,q_2\lambda_+,\alpha):=\frac{1}{\sqrt{\delta}}C(N,\mu,q_2\lambda_+,\alpha)(1+\sqrt{q_2})
$$
where $\delta$ and $C(N,\mu,q_2,\lambda_+,\alpha)$ is as in Proposition \ref{step 3}.
\end{remark}
\begin{proof}
We define the following rescaled function
$$
w(y) := \frac{u(x_0+\rho y )}{\Big ( \sqrt{\frac{2q_2}{\delta}}\cdot \rho+  \Linfty{u}{B_\rho(x_0)} + \frac{\rho}{\delta}\LN{F}{B_\rho(x_0)}\Big ) },\;\;y\in B_1
$$
where $F=|f_+|+|f_-|$.  We can easily verify that 
\be\label{linfty}
\Linfty{w}{B_1}\le 1.
\ee
From Remark \ref{rescaling} we check that $w$ is an absolute minimizer of the following functional 
\be\label{Jtilde}
\tilde{J} :=\int_{B_1} \Big (\bar A(y,w)|\nabla w|^2 -\bar f(y,w)w+\bar Q(y)\lambda(w)\Big )\,dy
\ee
where $\bar A_{\pm}, \bar f_{\pm}$ and $\bar Q$ are defined as follows
\[
\begin{split}
\bar A_{\pm}(y)&:= A_{\pm}(x_0+\rho y),\\
\bar f_{\pm}(y)&: = \frac{\rho^2\cdot f_{\pm}(x_0+\rho y)}{\left ( \sqrt{\frac{2q_2}{\delta}}\cdot \rho+\Linfty{u}{B_{\rho }(x_0)}+\frac{\rho}{\delta}\LN{F}{B_{\rho}(x_0)} \right )},\\
\bar Q(y)&= \frac{\rho^2\cdot Q(x_0+\rho y)}{\left ( \sqrt{\frac{2q_2}{\delta}}\cdot \rho+\Linfty{u}{B_{\rho}(x_0)}+\frac{\rho}{\delta}\LN{F}{B_{\rho}(x_0)} \right )^2},\;y\in B_1.
\end{split}
\]
We observe that the pair ($\tilde{J}$, $w$) satisfies the assumptions of Proposition \ref{step 3}. Indeed, we see from \eqref{SMALLNESS} that 
\be\label{coeff1.1}
\Linfty{\bar A_{\pm} - \bar A_{\pm}(0)}{B_1}=\|A_{\pm}-A_{\pm}(x_0)\|_{L^{\infty}(B_\rho(x_0))} \le \delta.
\ee
Also,
\be\label{coeff2.1}
\LN{\bar f_{\pm}}{B_1} = \frac{\rho \cdot \LN{f_{\pm}}{B_{\rho}(x_0)}}{\left ( \sqrt{\frac{2q_2}{\delta}}\cdot \rho+\Linfty{u}{B_{\rho}(x_0)}+\frac{\rho}{\delta}\LN{F}{B_{\rho}(x_0)} \right )}\le \delta.
\ee
and for the term $\bar Q$, since $\Linfty{Q-Q(0)}{B_{\rho}(x_0)}\le 2q_2$ we have 
\be\label{coeff3.1}
\Linfty{\b{Q} - \b{Q}(0)}{B_1} = \frac{\rho^2\cdot \Linfty{{Q} - {Q}(x_0)}{B_{\rho}(x_0)}}{\left ( \sqrt{\frac{2q_2}{\delta}}\cdot \rho+\Linfty{u}{B_{\rho}(x_0)}+\frac{\rho}{\delta}\LN{f}{B_{s_0}(\rho)}  \right )^2}\le \delta.
\ee
Also we can easily see that 
\be
|\b{Q}| \le \frac{\delta}{2} <1\le q_2 \mbox{ in $B_1$. }
\ee
Therefore from Proposition \ref{step 3} we have 
\be\label{estonw}
\begin{split}
\sup_{x\in B_{r}}|w(x)|\le C(N,\mu, q_2,\lambda_+,\alpha) r^{\alpha} \;\; \forall \, r\le R_0\\
\end{split}
\ee
and on plugging in the definition of $w$ in \eqref{estonw} we obtain for all $r\le R_0$
\be\label{estonu}
\begin{split}
\sup_{z\in B_{\rho r}(x_0)}|u(z)|&=\sup_{y\in B_{r}}|u(x_0+\rho y)|\\
&\le {C(N,\mu, q_2,\lambda_+,\alpha)} \left (  \sqrt{\frac{2q_2}{\delta}}\cdot \rho+\Linfty{u}{B_{\rho}(x_0)} +\rho \LN{f}{B_\rho(x_0)} \right ) r^{\alpha}.
\end{split}
\ee
By replacing $\rho r$ by $s$ in \eqref{estonu}, we obtain $\forall s\le \rho R_0$
\[
\begin{split}
\sup_{B_{s}(x_0)}|u(x)|&\le \frac{C(N, \mu, q_2,\lambda_+,\alpha)}{\rho^{\alpha}} \left (  \sqrt{\frac{2q_2}{\delta}}\cdot \rho+\Linfty{u}{B_{\rho}(x_0)} +\rho \LN{f}{B_\rho(x_0)} \right )  s^{\alpha} \\
&\le \frac{C_2(N, \mu, q_2,\lambda_+,\alpha)}{\rho^{\alpha}}\Big (1+\sqrt{q_2} \Big )\left (  \rho+\Linfty{u}{B_{\rho}(x_0)} +\rho \LN{f}{B_\rho(x_0)} \right )  s^{\alpha}.\;\; \\
\end{split}
\]
\end{proof}

\begin{remark}\label{rem 3.4}

{We have obtained a local asymptotic $C^{0,\alpha}$ regularity estimates on absolute minimizers $u$ in the balls centred at the zero level sets. Now suppose $x_0\in \{u>0\} \cap B_1$. }Since we know that $u$ is a continuous function in $B_1$, therefore the set $\{u>0\}$ is an open set. 

Suppose $B_{\rho}(x_0)\subset \subset \{u>0\}\cap B_1$. For any $\vf\in C_c^{\infty}(B_{\rho}(x_0))$ we observe that 
\be\label{t small}
|t|< \frac{\inf_{B_{\rho}(x_0)}(u)}{4\Linfty{\vf}{B_{\rho}(x_0)}}\implies u+t\vf >\frac{3\inf_{B_{\rho}(x_0)}(u)}{4}>0.
\ee
Hence $B_{\rho}(x_0)\subset \subset \{u+t\vf >0\}$ provided $t\in \R$ is as in \eqref{t small}.

We know that $u$ is an absolute minimizer of $J(\cdot;B_1)$, so is of $J(\cdot;{B_{\rho}(x_0)})$ (c.f. Lemma \ref{local minimizer}). Therefore 
\[
J(u;B_{\rho}(x_0)) \le J(u+t\vf;B_{\rho}(x_0)).
\]
Since $B_{\rho}(x_0) \subset \subset \{u>0\}$ and $B_{\rho}(x_0) \subset \subset \{u+t\vf>0\}$ for any $t$ such that $|t|< \frac{\min_{B_{\rho}(x_0)}(u)}{4\Linfty{\vf}{B_{\rho}(x_0)}}$ we have 
\[
\int_{B_{\rho}(x_0)} \Big ( A_+(x) |\nabla u|^2 +f_+ u\Big  )\,dx \le \int_{B_{\rho}(x_0)} \Big ( A_+(x) |\nabla (u+t\vf)|^2 +f_+ (u+t\vf )\Big )\,dx.
\]
In other words, for any $\vf\in C_c^{\infty}(B_{\rho}(x_0))$ 
\[
\frac{d}{dt}\Big|_{t=0} \int_{B_{\rho}(x_0)} \Big ( A_+(x) |\nabla (u+t\vf)|^2 +f_+ (u+t\vf )\Big )\,dx=0.
\]
Elaborating the above identity
\[
\begin{split}
\frac{d}{dt}\Big|_{t=0}\Bigg [ \int_{B_{\rho(x_0)}} \Big ( A_+(x)|\nabla (u+t\vf)|^2 +f_+(u+t\vf)  \Big )\,dx\Bigg ] &= \int_{B_{\rho(x_0)}} \big ( 2A_+(x) \nabla u\cdot \nabla \vf +f_+ \vf \big )\,dx\\
&=0.
\end{split}
\]
Therefore, $u$ is a weak solution to the following Euler-Lagrange PDE
\be \label{PDE+}
\dive\left ( A_+(x) \nabla u \right ) = -\frac{1}{2}f_+(x) \;\mbox{ in $B_{\rho}(x_0)$}.
\ee
Furthermore whenever $B_{\rho}(x_0) \subset \subset \{u<0\}\cap B_1$, proceeding similarly we obtain
\be \label{PDE-}
\dive\left ( A_-(x) \nabla u \right ) = -\frac{1}{2}f_-(x) \;\mbox{ in $B_{\rho}(x_0)$}.
\ee
In summary, the following PDEs are true in weak sense
$$
\begin{cases}
\dive(A_+(x)\nabla u)=-\frac{1}{2}f_+\qquad \mbox{in $\{u>0\}\cap B_{1}$}\\
\dive(A_-(x)\nabla u)=-\frac{1}{2}f_-\qquad \mbox{in $\{u<0\}\cap B_1$}.
\end{cases}
$$
\begin{comment}
From Proposition \ref{step 0.4} we know that $u$ is locally $C^{0,\alpha}$ in $\left ( \{u>0\}\cup \{u<0\} \right )\cap B_{1/2}$.  However, the regularity estimates on $u$ may deteriorate as we move closer to $F(u)\cap B_1$ (c.f. Corollary \ref{step 0.45}). Therefore we cannot yet conclude that $u\in C^{0,\alpha}(B_{1/2})$. In order to prove it, we utilize information on how the $C^{0,\alpha}$ estimates obtained in Proposition \ref{step 0.4} deteriorate near the free boundary (the Schauder estimates in Corollary \ref{step 0.45}). We also make use of the non-homogenous Moser-Harnack inequality along with some geometric arguments.
\end{comment}
\end{remark}

\begin{proposition}\label{step 4}
Let $u\in H^{1}(B_1)$ be a bounded absolute minimizer of $J(\cdot ;B_1)$.  Then for every $0<\alpha<1$ and for $\delta(N, \mu, q_2,\lambda_+,\alpha)>0$ as in Proposition \ref{step 3}, we have 
\be\label{claim 4}
\|A_{\pm}-A_{\pm}(0)\|_{L^{\infty}(B_1)} \le \frac{\delta}{2}  \implies \|u\|_{C^{\alpha}(B_{1/2})} \le C\left (1+\Linfty{u}{B_1}+\LN{F}{B_{1}}\right )
\ee
where $C:=C(N, \mu, q_2,\lambda_+,\alpha)$.
\end{proposition}
\begin{proof}
In this proof,  $R_0$ is as in Proposition \ref{step 3.5}. For the sake of this proof we introduce the following function in the set $\{u\neq 0\}\cap B_1$
$$
d(x) = \begin{cases}
\dist (x, \overline{ \{u\le 0\}})\;\;\mbox{if $u(x)>0$}\\
\dist (x,\overline{ \{u\ge 0\}})\;\;\mbox{if $u(x)<0$}.
\end{cases}
$$
We start by proving the following auxiliary estimates
\begin{enumerate}[label=\textbf{(e-\arabic*)}]
\item \label{e1} For any $y\in {B_{1/2}}$ and $x\in \{u=0\}\cap{ B_{5/8}}$ we have
\be\label{e1eq}
|u(x)-u(y)| = |u(y)| \le C_1(N,\mu,\alpha,q_2,\lambda_+ ) \Big ( 1+\Linfty{u}{B_1}+\LN{F}{B_1} \Big )|x-y|^{\alpha}.
\ee
\item \label{e1.5}
For any $x\in {B_{1/2}}$ 
\be\label{e1.5eq}
|u(x)| \le C_1(N,\mu,\alpha,q_2,\lambda_+ )\Big ( 1+\Linfty{u}{B_1}+\LN{F}{B_1} \Big ) d(x)^{\alpha}.
\ee
\item \label{e2} For any $x\in  \Big ( \{u>0\} \cup \{u<0\} \Big ) \cap B_{1/2}$ such that $d=d(x)\le \frac{R_0}{8}$ 
\be\label{e2eq}
\|u\|_{C^{0,\alpha}(B_{d/8}(x))} \le \frac{C_2(N,\mu,\alpha)}{d^{\alpha}} \Big ( u(x) + d \LN{F}{B_1} \Big ).
\ee
\item \label{e3} For any $x\in \Big ( \{u>0\} \cup \{u<0\} \Big ) \cap B_{1/2}$ such that $d=d(x)\le \frac{R_0}{8}$
\be\label{e3eq}
\|u\|_{C^{0,\alpha}(B_{d/8}(x))} \le C_3(N,\mu, q_2,\lambda_+,\alpha) \Big ( 1+ \Linfty{u}{B_1}+\LN{F}{B_1} \Big ).
\ee
\end{enumerate}
Before delving into the proofs of \ref{e1}-\ref{e3}. We start by observing that for any $x\in B_{1}$ we have
\be\label{smallness around x0}
\Linfty{A_{\pm}-A_{\pm}(x)}{B_{d}(x)} \le \Linfty{A_{\pm}-A_{\pm}(0)}{B_{d}(x)} + |A_{\pm}(0)-A_{\pm}(x)| \le \frac{\delta}{2} +\frac{\delta}{2} =\delta.
\ee
In order to prove \ref{e1}, we observe that $B_{1/4}(x)\subset \subset B_1$. Once $u$ is an absolute minimizer of $J(\cdot;B_1)$, so it is of $J(\cdot;B_{1/4}(x))$. We now divide the proof of \ref{e1} in two cases
\\[3pt]
\textbf {Case e1.A} : $y\in B_{1/2}$, $x\in B_{5/8}$ and $|x-y|< \frac{R_0}{4}$. 

\vspace{0.3cm}

In this case, since \eqref{smallness around x0} holds, the choice of $\rho=\frac{1}{4}$, $r=|x-y|$ and $x_0=x$ is an admissible choice in Proposition \ref{step 3.5}. It readily follows that for some constant $C_0:= C_0(N,\mu,q_2,\lambda_+,\alpha)$ 
$$
|u(x)-u(y)| = |u(y)| \le C_0(N,\mu,q_2,\lambda_+,\alpha ) \Big ( 1+\Linfty{u}{B_1}+\LN{F}{B_1} \Big )|x-y|^{\alpha}.
$$
\textbf{Case e1.B}: $y\in B_{1/2}$, $x\in B_{5/8}$ and $|x-y|\ge \frac{R_0}{4}$. 

$$
{|u(x)-u(y)|} \le \frac{4^\alpha\cdot 2\Linfty{u}{B_1}}{R_0^{\alpha}}|x-y|^{\alpha}.
$$
\noindent
Now, from the \textbf{Case e1.A} and \textbf{Case e1.B}, for $y\in B_{1/2}$ and $x\in \{u=0\}\cap B_{5/8}$ we have
$$
|u(x)-u(y)| \le \max \Big ( C_1, \frac{2\cdot 4^\alpha}{R_0^\alpha} \Big ) \Big ( 1+\Linfty{u}{B_1} + \LN{F}{B_1}\Big )|x-y|^{\alpha}.
$$
This proves \ref{e1} with $C_1 = \max \big ( C_0, \frac{2\cdot 4^\alpha}{R_0^\alpha}\big )$. 

In order to prove \ref{e1.5}, let us take $\bar x \in \{u=0\}$ be such that  $d=d(x) = |x-\bar x|$. We again divide the proof in two cases
\\[5pt]
\textbf{Case e2.A}: Assume $d(x)\le \frac{R_0}{8}$.

We observe that $\bar x\in B_{5/8}\cap \{u=0\}$, indeed,
$$
|\bar x| \le |x| + |x-\bar x| \le \frac{1}{2}+\frac{R_0}{8} < \frac{5}{8}.
$$
Also, we can easily verify that $B_{1/4}(\bar x) \subset \subset B_1$.  Hence using \ref{e1} with $\bar  x \in \{u=0\}\cap B_{5/8}$, we have 
\[
\begin{split}
|u(x)-u(\bar x)| = |u(x)| &\le C_{00}(N,\mu,q_2,\lambda_+,\alpha ) \Big ( 1+\Linfty{u}{B_1}+\LN{F}{B_1} \Big )|x-\bar x|^{\alpha}\\
&= C_{00}(N,\mu,q_2,\lambda_+,\alpha ) \Big ( 1+\Linfty{u}{B_1}+\LN{F}{B_1} \Big )d(x)^{\alpha}.
\end{split}
\]
\\
\textbf{Case e2.B}: Assume that $d>\frac{R_0}{8}$. In this case 
$$
|u(x)| \le \frac{8^\alpha \cdot \Linfty{u}{B_1}}{R_0^\alpha}d(x)^{\alpha}.
$$
\noindent
Thus, \textbf{Case e2.A} and \textbf{Case e2.B} prove \ref{e1.5} with $C_2 = \max \Big ( C_{00},  \frac{4^\alpha }{R_0^\alpha}\Big )$. 
\\[5pt]
For the proof of \ref{e2}, it is enough to consider only the case where $x\in B_{1/2}\cap \{u>0\}$, since the case where $x\in B_{1/2}\cap \{u<0\}$ can be treated similarly. From \eqref{PDE+} in Remark \ref{rem 3.4}, $u$ is a weak solution of the following PDE
\be\label{PDE++}
\dive\Big ( A_+(x)\nabla u \Big ) = \frac{1}{2}f_+ \;\;\mbox{in $B_{d/4}(x)$}.
\ee
From the  non-homogenous Moser-Harnack inequality \cite[Theorem 1]{serrin63}, we have  
\be\label{harnack}
\begin{split}
\sup_{B_{d/4}(x)} u& \le C (N,\mu) \left (  \inf_{B_{d/8}(x)}u+d\|f_+\|_{L^N(B_{d/4}(x))} \right )\\
&\le C (N,\mu) \left (  \inf_{B_{d/8}(x)}u+d\|F\|_{L^N(B_{d/4}(x))} \right ).
\end{split}
\ee
Moreover from \eqref{smallness around x0} and \eqref{estg2} in Corollary \ref{step 0.45} with $B_{\rho}(x_0) = B_{d/4}(x)$ and recalling that $u\ge 0$ in $B_{d/4}(x)$ we have
\be\label{schauder and harnack}
\begin{split}
\|u\|_{C^{0,\alpha}(B_{d/8}(x))} &\le \frac{C(N,\mu,\alpha)}{d^{\alpha}}  \left ( \sup_{B_{d/4}}(u)+d\LN{F}{B_{d/4}(x)} \right ).
\end{split}
\ee
Using \eqref{harnack} and \eqref{schauder and harnack} we arrive at
\[
\begin{split}
\|u\|_{C^{0,\alpha}(B_{d/8}(x))} &\le \frac{C (N,\mu,\alpha)}{d^{\alpha}} \left (  \inf_{B_{d/8}(x)}u+d\|F\|_{L^N(B_{d/4}(x))} \right )\\
&\le \frac{C (N,\mu,\alpha)}{d^{\alpha}} \left (  u(x)+d\|F\|_{L^N(B_{d/4}(x))} \right ).
\end{split}
\]
This concludes the proof of \ref{e2}. In order to prove \ref{e3}, we again treat only the case $x\in \{u>0\}\cap B_{1/2}$. From \ref{e1.5}
$$
u(x)=|u(x)| \le C(N,\mu,q_2,\lambda_+ ,\alpha)\Big ( 1+\Linfty{u}{B_1}+\LN{F}{B_1} \Big ) d(x)^{\alpha}.
$$
Plugging the above estimates in \ref{e2} we obtain 
\be\label{+phase}
\begin{split}
\|u\|_{C^{\alpha}(B_{d(x)/8}(x))}&\le C_1C_2 \left (1+\Linfty{u}{B_1}+\LN{F}{B_{1}} \right )+C_1d(x)^{1-\alpha}\|F\|_{L^N(B_1)}\\
&\le C_3(N, \mu, q_2,\lambda_+, \alpha) \left (1+\Linfty{u}{B_1}+\LN{F}{B_{1}} \right ).
\end{split}
\ee
Now, under the possession of \ref{e1}-\ref{e3}, we now finish the proof of Proposition \ref{step 4}. 

We again divide the proof in cases.
\\[5pt]
\textbf{Case I}.  Let $x,y\in B_{1/2}$ are such that $u(x)\cdot u(y)=0$. 
\\[5pt]
We can assume without loosing generality that $u(x)=0$. Then, it follows readily from \ref{e1} that 
$$
|u(x)-u(y)| = |u(y)| \le C_1(N,\mu,q_2,\lambda_+,\alpha ) \Big ( 1+\Linfty{u}{B_1}+\LN{F}{B_1} \Big )|x-y|^{\alpha}.
$$
\textbf{Case II}. Let $x,y\in B_{1/2}$ such that $u(x)\cdot u(y)\neq 0$. 
\\[5pt]
Without loss of generality, we can assume the 
$$
d(y)\le d(x).
$$
Once more, splitting the proof in cases
 
\vspace{0.2cm}
\noindent
\textbf{Case II.1}. If $|x-y|< \frac{d(x)}{8}$.  We now study the two subcases

\vspace{0.3cm}
\hspace{-0.5cm}
\textbf{Case II.1.A}. If $d(x)\le \frac{R_0}{8}$. 

\indent\indent
\hspace{-0.7cm}{In this case, $y\in B_{d(x)/8}(x)$.} 			Then it readily follows from \ref{e3} that 
			\[
			\begin{split}
			|u(x)-u(y)| &\le [u]_{C^{\alpha}(B_{d(x)/8}(x))}|x-y|^{\alpha}
			\le C_3 \Big ( 1+ \Linfty{u}{B_1}+\LN{F}{B_1} \Big )|x-y|^{\alpha}
			\end{split}
			\]
\indent \indent where $C_3:=C_3(N,\mu q_2,\lambda_+,\alpha,)$.	
\vspace{0.3cm}
		
\hspace{-0.5cm} \textbf{Case II.1.B}. If $d(x) > \frac{R_0}{8}$. 

\indent\indent
Since $\dive\big ( A_+(x)\nabla u \big ) = \frac{1}{2}f_+ \;\;\mbox{in $B_{d/8}(x)$}$ in weak sense, we apply \eqref{estg2} in Corollary \ref{step 0.45} 
\indent \indent on $u$ in the ball $B_{d/8}(x)$. This leads to 
		
			\[
			\begin{split}
			|u(x)-u(y)| &\le [u]_{C^{\alpha}(B_{d(x)/8}(x))}|x-y|^{\alpha}\\
			&\le \frac{C(N,\mu,\alpha)}{d^{\alpha}} \Big ( \Linfty{u}{B_1}+d \LN{F}{B_1} \Big )|x-y|^{\alpha}\\
			&\le \frac{C_4(N,\mu,\alpha)}{R_0^{\alpha}} \Big ( 1+ \Linfty{u}{B_1}+\LN{F}{B_1} \Big )|x-y|^{\alpha}.
			\end{split}
			\]

\textbf{Case II.2}. If $|x-y|\ge  \frac{d(x)}{8}$. 
\\[5pt]
\indent \indent By \ref{e1.5} and the assumption $d(x)\ge d(y)$ we obtain 
\be\label{more}
\begin{split}
|u(x)-u(y)| &= |u(x)|+|u(y)|\\
&\le  C_1(N, \mu, q_2,\lambda_+,\alpha) \Big (1+\Linfty{u}{B_1}+\LN{F}{B_{1}}\Big ) (d(x)^{\alpha}+d(y)^{\alpha})\\
&\le 2  C_1(N, \mu, q_2,\lambda_+,\alpha) \Big (1+\Linfty{u}{B_1}+\LN{F}{B_{1}}\Big )d(x)^{\alpha}\\
&\le C_5(N, \mu, q_2,\lambda_+,\alpha) \Big (1+\Linfty{u}{B_1}+\LN{F}{B_{1}}\Big )|x-y|^{\alpha}.
\end{split}
\ee
\noindent This proves Proposition \ref{step 4}.
\end{proof}
Now we present a scaled version of Proposition \ref{step 4}.

\begin{corollary}\label{step 4.5}
Let $u\in H^{1}(B_\rho(x_0))$ be a bounded absolute minimizer of $J(\cdot ;B_\rho(x_0))$ and $\rho< 1$.  Then for every $0<\alpha<1$ and for $\delta(N, \mu, q_2,\lambda_+,\alpha)>0$ as in Proposition \ref{step 3}, we have 
\be\label{claim 4}
\|A_{\pm}-A_{\pm}(x_0)\|_{L^{\infty}(B_\rho(x_0))} \le \frac{\delta}{2}  \implies \|u\|_{C^{\alpha}(B_{\frac{\rho}{2}}(x_0))} \le \frac{C}{\rho^{\alpha}}\left (\rho+\Linfty{u}{B_\rho(x_0)}+\rho \LN{F}{B_\rho(x_0)}\right )
\ee
where $C:=C(N, \mu, q_2,\lambda_+, \alpha)$.
\end{corollary}
\begin{proof}
We reduce the Corollary \ref{step 4.5} to Proposition \ref{step 4} by using  the following rescaling 
$$
w(y) = \frac{1}{\rho}u(x_0+\rho y).
$$
From Remark \ref{rescaling}, $w$ is an absolute minimizer of $\bar J(\cdot;B_1)$ given by 
$$
\bar J(v;B_1) = \int_{B_1} \Big ( \bar A(x,v)|\nabla v|^2 + \bar f(x,v) v + \bar Q(x)\lambda(v) \ \Big )\,dx
$$
where 
\[
\begin{split}
\bar A_{\pm}(y) &= A_{\pm}(x_0 + \rho y)\\
\bar f_{\pm}(y) &= \rho f_{\pm}(x_0+\rho y)\\
\bar Q(y)&= Q(x_0+\rho y).
\end{split}
\]
We observe that $w$ satisfies assumptions of the Proposition \ref{step 4}.  The final estimates are obtained by rescaling back to $u$.
\end{proof}
And now, we present the proof of the Theorem \ref{C01-}.
\section{Proof of the main result (Theorem \ref{C01-}).}\label{main theorem}

\begin{proof}[Proof of Theorem \ref{C01-}]
Let $D\subset \subset B_1$. We observe that $D\subset \subset B_{d} \subset \subset B_1$ where $d=1-\frac{\dist(D,\partial B_1)}{2}$. Since $A_{\pm}\in C(B_1)$, then $A_{\pm}$ are uniformly continuous in $\overline{B_d}$. Thus, we define $\o_{A_{\pm}, B_d}$ as the following  modulus of continuity (c.f. Definition \ref{mod of cont}). We define $\o_{A_{\pm},\overline{B_d}}$ to be 
$$
\o_{A_{\pm},\overline{B_d}}(t) := \max \Bigg ( \sup_{\substack{ |x-y|<t\\x,y\in \overline{B_d}}}|A_+(x)-A_+(y)|\,,\,\sup_{\substack{ |x-y|<t\\x,y\in \overline{B_d}}}|A_-(x)-A_-(y)|\Bigg )\mbox{  for $t\le 2d$}
$$
and we define $\o_{A_{\pm},\overline{B_d}}(t) := \o_{A_{\pm},\overline{B_d}}(2d)$ for $t>2d$.
{We set $t_0$ as  
$$
t_0 := t_0(\o_{A_{\pm},\overline{B_d}},\delta) = \sup \Big \{ t \;\Big | \;\o_{A_{\pm},\overline{B_d}}(t)\le \delta \Big \}
$$ 
as well as 
$$
s_0 := \min\Big (t_0,\frac{\dist(D,\partial B_1)}{4} \Big ).
$$
Since $\o_{A_{\pm},\overline {B_d}}$ is a non-decreasing function we have $\o_{A_{\pm},\overline{B_d}}(s_0)\le \delta$. }Now since 
$$
D\subset \bigcup_{x\in D}B_{s_0}(x) \subset \overline {B_d}
$$
we have 
\be\label{smallness1...}
\sup_{B_{s_0}(x)}|A_{\pm}-A_{\pm}(x_0)| \le \o_{A_{\pm},\overline{B_d}}(s_0) \le {\delta},\;\forall x\in D.
\ee
By Lemma \ref{local minimizer} $u$ is also an absolute minimizer of $J(\cdot;B_{s_0}(x))$, we have from Corollary \ref{step 4.5} that for all $y\in B_{s_0/2}(x)\cap D$ 
\be\label{near x_0...}
|u(x)-u(y)|\le \frac{C(N, \mu,q_2, \lambda_+,\alpha,)}{s_0^{\alpha}}\left (s_0+\Linfty{u}{B_{1}}+s_0\LN{f}{B_{1}}  \right )|x-y|^{\alpha}.
\ee
Now, if $x, y\in D$ are such that $|x-y| \ge s_0/2$, then 
\be\label{far from x_0...}
{|u(x)-u(y)|}\le 2^{1+\alpha} \frac{\Linfty{u}{B_1}}{s_0^{\alpha}}|x-y|^{\alpha}.
\ee
By combining \eqref{near x_0...} and \eqref{far from x_0...}, we arrive to $\Big ($since $s_0\le\frac{\dist(D,\partial B_1)}{4}\le\frac{\diam(B_1)}{4} =\frac{1}{2}<1$$\Big )$
\be\label{precise est...}
[u]_{C^{\alpha}(D)} \le \frac{C(N,\alpha, \mu, q_2,\lambda_+)}{s_0^{\alpha}}\left ( 1+\Linfty{u}{B_{1}}+\LN{f}{B_{1}} \right ).
\ee
{We observe from the definition of $s_0$ that 
$$
[u]_{C^{\alpha}(D)} \le \begin{cases}
\frac{C(N, \mu,q_2,\lambda_+,\alpha)}{t_0^{\alpha}}\left (1+ \Linfty{u}{B_{1}}+\LN{f}{B_{1}} \right ) \;\;\;\;\;\;\;\;\;\;\;\;\;\;\;\;\;\;\;\mbox{if $t_0\le \frac{\dist(D,\partial B_1)}{4}$}\\[7pt]
\frac{4^\alpha\cdot C(N, \mu,q_2,\lambda_+,\alpha)}{\dist(D,\partial B_1)^\alpha}\left ( 1+\Linfty{u}{B_{1}}+\LN{f}{B_{1}} \right )\;\; \;\;\;\;\;\;\;\;\;\;\;\;\;\mbox{if $t_0\ge \frac{\dist(D,\partial B_1)}{4}$.}
\end{cases}
$$
}
{In order to control the first term in the equation above by a universal multiple of $\dist(D,\partial B_1)^{-\alpha}$, we observe that} once $t_0>0$ depends only on the modulus of continuity $\o_{A_{\pm}, \overline{B_d}}$ and $\delta$, there exists a $n_0:= n_0(\o_{A_{\pm}, \overline{B_d}},\delta)=n_0(N,\mu,\alpha,q_2,\lambda_+,\o_{A_{\pm},\overline{B_d}})\in \N$ such that  $\frac{2}{n_0}\le t_0$. Hence
$$
\frac{\dist(D,\partial B_1)}{n_0}\le \frac{\diam(B_1)}{n_0} =\frac{2}{n_0}\le t_0 \implies \frac{1}{t_0^{\alpha}}\le \frac{n_0^{\alpha} }{\dist(D,\partial B_1)^{\alpha}}= \frac{C(N,\mu,q_2,\lambda_+,\alpha,\o_{A_{\pm},\overline{B_d}}) }{\dist(D,\partial B_1)^{\alpha}}.
$$
\begin{comment}
Thus the estimate \eqref{precise est} becomes
\be\label{last est}
[u]_{C^{\alpha}(D)} \le \frac{C(N, \alpha, \mu, q_2\lambda_+, \o_{A_{\pm}, D})} {\dist(D,\partial B_1)^{\alpha}}\left ( 1+\Linfty{u}{B_1}+\LN{f}{B_1} \right ).
\ee
We consider a finite cover of $\overline D\subset \subset B_1$ 
$$
\overline D \subset \bigcup_{i=1}^n \Big \{ B_{\frac{s_0}{8}}(x_i): \,x_i\in D \Big \}.
$$
Let $x,y\in D$, and assume $x\in B_{\frac{s_0}{8}(x_i)}$ for some $i\in \{1,...,n\}$. If $|x-y|\le \frac{s_0}{16}$, then 
$$
|y-x_i|\le |x-y|+|x-x_i|\le \frac{s_0}{16}+\frac{s_0}{8}= \frac{3s_0}{19}<<\frac{s_0}{2}.
$$
Therefore, $x,y\in B_{s_0/2}(x_i)$. From \eqref{last est} with $x_0=x_i$ we obtain
\be\label{x,y near}
|u(x)-u(y)| \le \frac{C(N, \alpha, \mu, q_2\lambda_+, \o_{A_{\pm}, D})} {\dist(D,\partial B_1)^{\alpha}}\left ( 1+\Linfty{u}{B_1}+\LN{f}{B_1} \right )|x-y|^{\alpha}
\ee
If $|x-y| \ge \frac{s_0}{8}$ then from \eqref{s0} we have
\be\label{x,y far}
|u(x)-u(y)| \le \frac{8^{\alpha}\cdot 2\Linfty{u}{B_1}}{s_0^{\alpha}}|x-y|^{\alpha}\le \frac{C(\o_{A_{\pm},D},\alpha)\Linfty{u}{B_1}}{\dist(D,\partial B_1)^{\alpha}}|x-y|^{\alpha}.
\ee
\end{comment}
Now, \eqref{precise est...} becomes
\be\label{general est...}
[u]_{C^{\alpha}(D)} \le \frac{C(N,\mu, q_2,\lambda_+,\alpha,\o_{A_{\pm},\overline{B_d}})}{\dist(D,\partial B_1)^{\alpha}}\left (1+ \Linfty{u}{B_{1}}+\LN{f}{B_{1}} \right ).
\ee
To finish we observe that Theorem \ref{C01-} follows from \eqref{general est...} by taking $D$ as $B_r$ for any $r<1$, once $\dist(D,\partial B_1)^{\alpha} = \dist(B_r,\partial B_1)^{\alpha} = (1-r)^{\alpha}$ and $d=\frac{1+r}{2}$. 

\end{proof}

 \newcommand{\noop}[1]{}

\ack
The work of the first author is supported by CNPq-311566/2019-7 (Brazil).

\end{document}